\documentclass[12pt]{elsarticle}
\usepackage{subfig, caption, keyval}
\usepackage{multirow}
\usepackage{color}
\usepackage{latexsym,amsmath,amsfonts,amscd,amsthm}
\usepackage{epsfig}
\usepackage{changebar}
\numberwithin{equation}{section}

\theoremstyle{plain}
\newtheorem{thm}{Theorem}[section]
\newtheorem{lem}[thm]{Lemma}

\newtheorem{prop}[thm]{Proposition}

\newtheorem{de}[thm]{Definition}

\theoremstyle{definition}

\newcommand{\BigOh}{\mathcal{O}}

\begin{document}
\baselineskip=1.5pc

\vspace{.5in}

\begin{center}

{\large\bf High order operator splitting methods based on an integral deferred correction framework }

\end{center}

\vspace{.1in}

\centerline{
Andrew J. Christlieb \footnote{Department of Mathematics and Department of Electrical and Computer Engineering,
Michigan State University, East Lansing, MI 48824, USA.
E-mail: christli@msu.edu} \qquad
Yuan Liu\footnote{Department of Mathematics,
Michigan State University, East Lansing, MI 48824, USA.
E-mail: yliu7@math.msu.edu} \qquad
 Zhengfu Xu\footnote{Department of Mathematical Science, Michigan Technological University, Houghton, MI 49931, USA. E-mail: zhengfux@mtu.edu.
}
}

\vspace{.6in}

\centerline{\bf Abstract}

\vspace{.1in}
Integral deferred correction (IDC) methods have been shown to be an efficient way to achieve arbitrary high order accuracy and possess good stability properties. In this paper, we construct  high order operator splitting schemes using the  IDC procedure to solve initial value problems (IVPs). We present analysis to show that the IDC methods can correct for both the splitting  and numerical errors, lifting the order of accuracy by $r$ with each correction, where $r$ is the order of accuracy of the method used to solve the correction equation. We further apply this framework to solve partial differential equations (PDEs). Numerical examples in two dimensions of linear and nonlinear initial-boundary value problems are presented to demonstrate the performance of the proposed IDC approach.

\vfill

{\bf Key Words:}
Integral deferred correction, initial-boundary value problem,
high-order accuracy, operator splitting.

%{\bf AMS(MOS) subject classification:} 65M99

\pagenumbering{arabic}

%section 1
%\input{introduction.tex}

\section{Introduction}
\label{sec1}
\setcounter{equation}{0}
\setcounter{figure}{0}
\setcounter{table}{0}

In this paper we present high order operator splitting methods based on the integral deferred correction (IDC) mechanism.
% The methods are motivated by a host of
%interesting current issues in scientific computing including: efficient multi-core implicit solvers \cite{Christlieb_Ong};  addressing  complex geometry on a
%uniform mesh \cite{MOLT1,MOLT2}; and combining parabolic and hyperbolic solvers in a high order adaptive mesh refinement (AMR) framework \cite{Claridge,shen2011adaptive,gibou2013high}.
The methods are designed to leverage recent progress on parallel time stepping and offer a great
deal of flexibility for computing the ordinary differential equations (ODEs). We focus on extending IDC theory to the case of splitting schemes on the IVP
\begin{equation}
\label{eqn:split-ode}
u_{t} = f(t,u) = \sum^{\Lambda}_{\nu = 1} f_{\nu}(t,u),  \qquad u(0) = u_0, \qquad t\in [0,T],
\end{equation}
and discuss the application in parabolic PDEs. Here, $u \in \mathbb{R}^n$ and $f(t, u) : \mathbb{R}^{+}\times \mathbb{R}^n \rightarrow \mathbb{R}^n$.

In the case that (\ref{eqn:split-ode}) 
%is a stiff ODE system, explicit numerical schemes suffer a severe restriction on the time step size. To overcome this difficulty,
%implicit time stepping methods are often used to allow for much larger time steps. Moreover, when 
arises from a method of lines discretization of time dependent PDEs which describe multi-physics problems, we encounter high dimensional computation. For these problems, the splitting methods can be applied to decouple the problems into simpler sub-problems. Therefore, the main advantage of operator splitting methods are problem simplification, dimension reduction,
and lower computational cost. Two broad categories can classify many splitting methods: differential operator splitting \cite{Bagrinovskii, Marchuk, Strang1, Strang2}
and algebraic splitting with the prominent example of the alternating direction implicit(ADI) method, which was first introduced in \cite{Doug2, Doug1,Peachman} for solving two dimensional heat equations. The main barrier in designing high order numerical methods based on the idea of splitting is the operator splitting error. To obtain high order accuracy via low order splitting method generally adds complexity to designing a scheme and stability analysis \cite{Hansen, McLachlan, Jia, McLachlan, McLachlan1, THALHAMMER, Geiser}.
%
%% ------------------------------------------------------------------------- %%
%%
%% ------------------------------------------------------------------------- %%
 A recent work in \cite{Bourlioux} utilizes the spectral deferred correction (SDC) procedure to the
advection-diffusion-reaction system in one dimension in order to enhance the overall order of accuracy.
However, their work does not contain a proof that the corrections raise the order
of the method.
%% ------------------------------------------------------------------------- %%

%% ------------------------------------------------------------------------- %%
%%
%% ------------------------------------------------------------------------- %%

In \cite{Dutt}, a SDC method is first proposed as a new variation on the classical deferred correction methods \cite{Bohmer}. The key idea is to recast the error equation such that the residual appears in the error equation in  integral form instead of differential form, which greatly stabilizes the method. It is proposed as a framework to generate arbitrarily high order methods. This family of methods use Gaussian quadrature nodes in the correction to the defect or error, hence the method can achieve a maximal order of $2(M-1)$ on $M$ grid points with $2(M-2)$ corrections. This main feature of the SDC method made it popular and extensive investigation can be found in \cite{Dutt,Minion,Layton3, Layton4,Layton1,Huang1,
Huang2, liu2008strong}. Following this line of approach, the IDC methods are introduced in \cite{Andrew, Andrew1,Andrew2}.  High order explicit and implicit Runge-Kutta (RK) integrators in both the prediction and correction steps (IDC-RK) are developed by utilizing uniform
quadrature nodes for computing the residual. In \cite{Andrew, Andrew1}, it is
established that using explicit RK methods of order $r$ in the
correction step results in $r$ higher degrees of accuracy with each successive
correction step, but only if uniform nodes are used instead of the Gaussian
quadrature nodes of SDC. It is shown in \cite{Andrew1} that the new methods
 produced by the  IDC procedure are yet again RK methods.
It is also demonstrated that, for the same order, IDC-RK methods possess better
stability properties than the equivalent SDC methods.  Furthermore, for explicit
methods, each correction of IDC or SDC increases the region of absolute
stability.  Similar results are generalized to  arbitrary order implicit and
additive RK methods in \cite{Andrew2}. Generally, for \emph{implicit}
methods based on IDC and SDC, the stability region becomes smaller when more
correction steps are employed.  It is believed that this is due to the numerical
approximation of the residual integral. The primary purpose of this work is to apply the IDC
methods to the low order operator splitting methods in order to obtain higher order accuracy.

%% ------------------------------------------------------------------------- %%
The paper is organized as follows. In Section 2, we briefly review several
classical operator splitting methods and show how these methods can be cast as
additive RK (ARK) methods.  In Section 3 we formulate the IDC
methodology for application to operator splitting schemes.  In Section 4, we prove that
IDC methods can correct for both the splitting and
numerical errors of ODEs, giving $r$ higher degrees of accuracy with each correction,
where $r$ is the order of the method used in the correction steps.  In section 5, as an interesting example, we will show how to use integral deferred correction for operator splitting (IDC-OS) schemes as a temporal discretization when solving PDEs.  In Section
6 we carry out numerical simulations  based on IDC methods for
both linear and non-linear parabolic equations, and demonstrate that the new
framework can achieve high order accuracy in time.  In Section 7 we conclude the
paper and discuss future work.  We note that both the parallel time stepping
version of IDC and the work presented in this paper are likely to benefit
from the work in \cite{Rokhlin}, and will be the subject of further investigation.
%% ------------------------------------------------------------------------- %% 

% section 2
%\input{section2.tex}
\section{Operator splitting schemes for ODEs}
\label{section2}

In this section, we review several splitting methods which will serve as the base solver in the IDC framework. For differential operator splitting, such as Lie-Trotter splitting and Strang splitting,  which happens at continuous level, we will apply appropriate numerical methods to the sub-problems and refer the whole approach as the discrete form of differential splitting. For both the differential splitting and algebraic splitting, we will show that each of the numerical schemes can be
written as an ARK method.  This insight is the first step
required to apply the IDC methodology \cite{Andrew2} to operator splitting schemes,
which is the primary purpose of the present work.
%% ------------------------------------------------------------------------- %%

\subsection{Review of ARK methods}

For IVP (\ref{eqn:split-ode}), when different $p$-stage RK integrators are applied to each operator $L_{\nu}$, the entire numerical method is called an ARK method. If we define the numerical solution after $n$ time steps as $\upsilon^n$, which is an approximation to the exact solution $u(t_n)$, then one step of a $p$-stage ARK method is given by
%% ------------------------------------------------------------------------- %%

%% ------------------------------------------------------------------------- %%
\begin{eqnarray}
\label{eqn:ark}
& \displaystyle{ \upsilon^{n+1} = \upsilon^n + \Delta t \sum^{\Lambda}_{\nu =1}\sum^p_{i=1}b_i^{[\nu]} f_{\nu}(t_n+c^{[\nu]}_i \Delta t, \tilde{\upsilon}_i)}, \\
\text{with}
 & \displaystyle{ \tilde{\upsilon}_i = \upsilon^n + \Delta t\sum^{\Lambda}_{\nu =1}\sum^p_{j = 1} a^{[\nu]}_{ij} f_{\nu}(t_n+c^{[\nu]}_j \Delta t, \tilde{\upsilon}_j)}.
\end{eqnarray}
%% ------------------------------------------------------------------------- %%
%
%% ------------------------------------------------------------------------- %%
and $\Delta t = t_{n+1}-t_n$. An ARK method is succinctly identified by its Butcher
tableau, as is demonstrated in Table \ref{table1}.
%
%% ------------------------------------------------------------------------- %%
%%
%% ------------------------------------------------------------------------- %%
\begin{table}[!htb]
\centering
\bigskip
\begin{tabular}{c| c| c | cccc | cc c| c ccc  }
 $c^{[1]}_1$   & $\cdots$ &  $c^{[\Lambda]}_1$  &  $a^{[1]}_{11}$ & $a^{[1]}_{12}$ & $\cdots$ & $a^{[1]}_{1p}$ &  &  $\cdots$ &   &  $a^{[\Lambda]}_{11}$ & $a^{[\Lambda]}_{12}$ & $\cdots$ & $a^{[\Lambda]}_{1p}$ \\
 $c^{[1]}_2$    & $\cdots$  &  $c^{[\Lambda]}_2$ &  $a^{[1]}_{21}$ & $a^{[1]}_{22}$ & $\cdots$ & $a^{[1]}_{2p}$ &  &  $\cdots$ &   &  $a^{[\Lambda]}_{21}$ & $a^{[\Lambda]}_{22}$ & $\cdots$ & $a^{[\Lambda]}_{2p}$ \\
 $\vdots$  & $\cdots$  &  $\vdots$  &  $\vdots$       & $\vdots$       & $\ddots$ & $\vdots $      &  &  $\cdots$ &   &  $\vdots$       & $\vdots$       & $\ddots$ & $\vdots$      \\
 $c^{[1]}_p$   & $\cdots$ & $c^{[\Lambda]}_p$ &  $a^{[1]}_{p1}$ & $a^{[1]}_{p2}$ & $\cdots$ & $a^{[1]}_{pp}$ &  &  $\cdots$ &   &  $a^{[\Lambda]}_{p1}$ & $a^{[\Lambda]}_{p2}$ & $\cdots$ & $a^{[\Lambda]}_{pp}$ \\ \hline
    \multicolumn{3}{c|}{}         &  $b^{[1]}_{1}$  & $b^{[1]}_{2}$  & $\cdots$ & $b^{[1]}_{p}$  &  & $\cdots$  &   &  $b^{[\Lambda]}_{1}$  & $b^{[\Lambda]}_{2}$  & $\cdots$ & $b^{[\Lambda]}_{p}$   \\
   \end{tabular}
   \caption{Butcher tableau for a $p$-stage ARK method. }
   \label{table1}
\end{table}

%% ------------------------------------------------------------------------- %%
In the following sections, we will explicitly write out the Butcher tableau
for each operator splitting scheme and conclude that each of the operator
splitting schemes considered in this work is indeed a form of ARK method.
%% ------------------------------------------------------------------------- %%

%% ------------------------------------------------------------------------- %%
%%
%% ------------------------------------------------------------------------- %%
\subsection{Lie-Trotter splitting}

%% ------------------------------------------------------------------------- %%
We describe Lie-Trotter splitting for \eqref{eqn:split-ode} in the case of
$\Lambda = 2$ in the right hand side functions.  We consider a single interval $[t_n, t_{n+1}]$ .  With first order Lie-Trotter splitting,
\eqref{eqn:split-ode} can be solved by two sub-problems:
\begin{equation}
   \label{3.2}
\left\{
  \begin{array}{ll}
  \vspace{0.05in}
     u_t = f_1(t, u), & \hbox{on  $ [t_n, t_{n+1}]$,} \\
     u_t = f_2(t, u), & \hbox{on  $[t_n, t_{n+1}]$.}
  \end{array}
\right.
\end{equation}
%% ------------------------------------------------------------------------- %%
The solution calculated from the first equation is used as
the initial value of the second equation. Note that this splitting occurs on the
continuous level. In order to define a discrete solver for \eqref{eqn:split-ode}, we need to
choose a numerical scheme for solving each sub-problem.  For example, if we use the
backward Euler scheme to solve both equations, we obtain a scheme
of the form
%% ------------------------------------------------------------------------- %%
\begin{equation}
   \label{3.3}
\left\{
  \begin{array}{ll}
  \vspace{0.05in}
    \displaystyle{ \frac{\widetilde{\upsilon}-\upsilon^n}{\Delta t}} = f_1(t_{n+1},\widetilde{\upsilon}), & \hbox{} \\
     \displaystyle{ \frac{\upsilon^{n+1}-\widetilde{\upsilon}}{\Delta t}}  = f_2(t_{n+1}, {\upsilon}^{n+1}) ~,~ & \hbox{}
  \end{array}
\right.
\end{equation}
where $\upsilon^n$ denotes the numerical approximation for $u$ at time level $t_n$. However, this approach only produces a first order approximation.
%% ------------------------------------------------------------------------- %%

%% ------------------------------------------------------------------------- %%
In order to make use of IDC methodology \cite{Andrew2} to lift the order of
accuracy of \eqref{3.3}, we write a Butcher tableau
for \eqref{3.3} in Table \ref{table2}.
Comparing the Butcher tableau for the Lie-Trotter splitting  with the general
form of the  Butcher tableau of an ARK method,  we can
view the discrete form of Lie-Trotter splitting (\ref{3.3})  as a 2-stage
ARK method.  This can be
extended to the case of $\Lambda$ operators, where the resulting Butcher tableau would
be a $\Lambda$-stage ARK method.
%% ------------------------------------------------------------------------- %%

%\begin{table}[!htb]
%\begin{center}
%\bigskip
%\begin{tabular}{c | cc | c c  }
% 1 &  1      &   0    &     0     &    0     \\
% 1 &  1      &   0    &     0     &    1      \\ \hline
%   &  1      &   0    &     0     &    1      \\
% %\hline
%   \end{tabular}
%   \label{table2}
%   \caption{Butcher tableau for Lie-Trotter splitting }
%\end{center}
%\end{table}

%% ------------------------------------------------------------------------- %%
\begin{table}
\begin{center}
\bigskip
\begin{tabular}{c | ccc | cc c  }
  0 & 0       &   0    &     0     &    0  & 0 & 0    \\
  1 & 0       &   1    &     0     &    0  & 0 & 0    \\
  1 & 0       &   1    &     0     &    0  & 0 & 1     \\ \hline
     &  0      &   1    &     0     &    0   & 0  & 1   \\
 %\hline
   \end{tabular}
      \caption{Butcher tableau for Lie-Trotter splitting.
      %This table is based
%      on first-order Lie-Trotter splitting, where each sub-problem is solved
%      using fully implicit backward Euler method.  We remark that it's possible
%      to define a Butcher tableau for operators with $\Lambda > 2$ right hand
%      side values, which would require $\Lambda$ different arrays.  Here, we
%      only present the case of $\Lambda = 2$.
}
   \label{table2}
\end{center}
\end{table}
%% ------------------------------------------------------------------------- %%

%In order to solve equation (\ref{3.3}) with Dirichlet boundary condition, we use
%\begin{equation}
%\label{3.4}
%\widetilde{U}^{n+1}(x,y) = g(t^{n+1},x,y)
%\end{equation}
%as boundary values when $(x,y)$ is on the boundary.

%% ------------------------------------------------------------------------- %%
%%
%% ------------------------------------------------------------------------- %%
\subsection{Strang splitting}

%% ------------------------------------------------------------------------- %%
In this section, we consider the second order Strang splitting for the case of
three operators to demonstrate how to construct Butcher tableaus for general differential
operator splitting schemes.
The case of $\Lambda = 3$ operators can arise when splitting a stiff ODE into three sub-problems while maintaining second order accuracy in time.

We also focus on a single time step, $[t_n, t_{n+1}]$.
Second order Strang splitting for \eqref{eqn:split-ode} reads as
\begin{align}
\begin{cases}\label{3.6}
 \vspace{0.1in}
 \displaystyle{u_t}=  f_1(t, u), \qquad  t \in [t_n,t_{n+\frac{1}{2}}], \\
                                 \vspace{0.1in}
\displaystyle{u_t}=  f_2(t, u), \qquad  t \in [t_{n+\frac{1}{2}},t_{n+1}], \\
                                 \vspace{0.1in}
\displaystyle{u_t}=  f_3(t, u), \qquad  t \in [t_n,t_{n+1}], \\
                                 \vspace{0.1in}
\displaystyle{u_t}=  f_2(t, u), \qquad  t \in [t_n,t_{n+\frac{1}{2}}], \\
                                 \vspace{0.1in}
\displaystyle{u_t}=  f_1(t, u), \qquad  t \in [t_{n+\frac{1}{2}}, t_{n+1}].
\end{cases}
\end{align}
Note that this splitting occurs on the continuous level, i.e. the temporal
derivative for each sub-problem in \eqref{3.6} has yet to be discretized.
If we discretize  equations (\ref{3.6})
with trapezoidal rule, we obtain an update of the form,
\begin{align}
\begin{cases}\label{3.2.1}
\vspace{0.1in}
\displaystyle{\frac{\tilde{\upsilon}_1 -\upsilon^n}{\frac{1}{2}\Delta t}= \frac{1}{2}( f_1(t_n,\upsilon^n)+ f_1(t_{n+\frac{1}{2}},\tilde{\upsilon}_1)) },  \\
                                 \vspace{0.1in}
\displaystyle{\frac{\tilde{\upsilon}_2-\tilde{\upsilon}_1}{\frac{1}{2}\Delta t}=  \frac{1}{2}(f_2(t_{n+\frac{1}{2}},\tilde{\upsilon}_1)+f_2(t_{n+1},\tilde{\upsilon}_2))},  \\
                                 \vspace{0.1in}
\displaystyle{\frac{\tilde{\upsilon}_3-\tilde{\upsilon}_2}{\Delta t}= \frac{1}{2}( f_3(t_n,\tilde{\upsilon}_2)+f_3(t_{n+1}, \tilde{\upsilon}_3))}, \\
                                 \vspace{0.1in}
\displaystyle{\frac{\tilde{\upsilon}_4-\tilde{\upsilon}_3}{\frac{1}{2}\Delta t}= \frac{1}{2}( f_2(t_n, \tilde{\upsilon}_3)+f_2(t_{n+\frac{1}{2}},\tilde{\upsilon}_4))} , \\
                                 \vspace{0.1in}
\displaystyle{\frac{\upsilon^{n+1}-\tilde{\upsilon}_4}{\frac{1}{2}\Delta t}= \frac{1}{2}( f_1(t_{n+\frac{1}{2}},\tilde{\upsilon}_4)+f_1(t_{n+1},\upsilon^{n+1}) )},
\end{cases}
\end{align}
where $t_{n+\frac{1}{2}} = t_n +\frac{1}{2}\Delta t$. In Table \ref{table3}, we write this scheme in the Butcher tableau.  Comparing this with the Butcher tableau of the ARK methods, again we see that we can view the Strang splitting as a 5-stage ARK scheme.
%% ------------------------------------------------------------------------- %%

\begin{table}[!htb]
\centering
\bigskip
%\begin{small}
\begin{tabular}{c| c| c | cccccc | ccc ccc| ccc ccc  }
 0   & 0              &  0  &  0        & 0  & 0  & 0  &  0 & 0 & 0 &  0 &  0 &  0 &  0 & 0  & 0 &0 &0 &0 &0 &0 \\[4pt]
  $\frac{1}{2}$   & $\frac{1}{2}$    &  0  &  $\frac{1}{4}$        & $\frac{1}{4}$  & 0  & 0  &  0 & 0 & 0 &  0 &  0 &  0 &  0 & 0  & 0 &0 &0 &0 &0 &0 \\[4pt]
  0   &1  &  0  & $\frac{1}{4}$  & $\frac{1}{4}$  & 0  & 0  &  0 & 0 & 0 &  $\frac{1}{4}$ &  $\frac{1}{4}$ &  0 &  0 & 0  & 0 &0 &0 &0 &0 &0 \\[4pt]
  0   &0  &  1  & $\frac{1}{4}$  & $\frac{1}{4}$  & 0  & 0  &  0 & 0 & 0 &  $\frac{1}{4}$ &  $\frac{1}{4}$ &  0 &  0 & 0  & 0 &0 &$\frac{1}{2}$ &$\frac{1}{2}$ &0 &0 \\[4pt]
  $\frac{1}{2}$   &$\frac{1}{2}$  &  0  & $\frac{1}{4}$  & $\frac{1}{4}$  & 0  & 0  &  0 & 0 & 0 &  $\frac{1}{4}$ &  $\frac{1}{4}$ &  $\frac{1}{4}$ &  $\frac{1}{4}$ & 0  & 0 &0 &$\frac{1}{2}$ &$\frac{1}{2}$ &0 &0 \\[4pt]
  1   &0 &  0  & $\frac{1}{4}$  & $\frac{1}{4}$  & 0  & 0  &  $\frac{1}{4}$ & $\frac{1}{4}$ & 0 &  $\frac{1}{4}$ &  $\frac{1}{4}$ &  $\frac{1}{4}$ &  $\frac{1}{4}$ & 0  & 0 &0 &$\frac{1}{2}$ &$\frac{1}{2}$ &0 &0 \\[4pt]
 \hline
    \multicolumn{3}{c|}{}  &  $\frac{1}{4}$  & $\frac{1}{4}$  & 0 & 0  &  $\frac{1}{4}$ & $\frac{1}{4}$ & 0 & $\frac{1}{4}$ & $\frac{1}{4}$  & $\frac{1}{4}$&  $\frac{1}{4}$  & 0  & 0 & 0 & $\frac{1}{2}$   & $\frac{1}{2}$ & 0 & 0 \\[4pt]
 %\hline
   \end{tabular}
   \caption{Butcher tableau for Strang splitting with $\Lambda = 3$.    }
   \label{table3}
%\end{small}
\end{table}
%% ------------------------------------------------------------------------- %%

%% TODO - I'm not sure you want to say this here ... -DS
%When we apply Strang splitting to the 2D equation (\ref{eqn:parabolic}), the third spatial operator $L_3$ in equation (\ref{3.6}) is assumed to be 0. Hence, the Butcher tableau will be much simpler, as shown in Table \ref{table4}. We note that the boundary values of the PDE are assigned according to the time stage via the appropriate Dirichlet boundary condition.
%The case of $\Lambda = 2$ operators arises when splitting \eqref{eqn:parabolic}
%by dimension.  In this case, the structure of the Butcher tableau simplifies by setting $L_3 = 0$,
%that is presented in Table \ref{table4}.
%% ------------------------------------------------------------------------- %%

%% ------------------------------------------------------------------------- %%
%% Table concerning 2 operators + Strang splitting.
%% ------------------------------------------------------------------------- %%
\begin{table}[!htb]
\centering
\bigskip
%\begin{small}
\begin{tabular}{c| c| cccccc | ccc ccc }
 0   & 0              &  0          & 0  & 0  & 0  &  0 & 0 & 0 &  0 &  0 &  0 &  0 & 0   \\[4pt]
  $\frac{1}{2}$   & $\frac{1}{2}$      &  $\frac{1}{4}$        & $\frac{1}{4}$  & 0  & 0  &  0 & 0 & 0 &  0 &  0 &  0 &  0 & 0   \\[4pt]
  0   &1    & $\frac{1}{4}$  & $\frac{1}{4}$  & 0  & 0  &  0 & 0 & 0 &  $\frac{1}{4}$ &  $\frac{1}{4}$ &  0 &  0 & 0  \\[4pt]
  $\frac{1}{2}$   &$\frac{1}{2}$    & $\frac{1}{4}$  & $\frac{1}{4}$  & 0  & 0  &  0 & 0 & 0 &  $\frac{1}{4}$ &  $\frac{1}{4}$ &  $\frac{1}{4}$ &  $\frac{1}{4}$ & 0   \\[4pt]
  1   &0   & $\frac{1}{4}$  & $\frac{1}{4}$  & 0  & 0  &  $\frac{1}{4}$ & $\frac{1}{4}$ & 0 &  $\frac{1}{4}$ &  $\frac{1}{4}$ &  $\frac{1}{4}$ &  $\frac{1}{4}$ & 0  \\[4pt]
 \hline
    \multicolumn{2}{c|}{}  &  $\frac{1}{4}$  & $\frac{1}{4}$   & 0  &  $\frac{1}{4}$ & $\frac{1}{4}$ & 0 & $\frac{1}{4}$ & $\frac{1}{4}$  & $\frac{1}{4}$&  $\frac{1}{4}$  & 0  &0 \\[4pt]
 %\hline
   \end{tabular}
   \caption{Butcher tableau for Strang splitting when $L_3 = 0$. }
   \label{table4}
%\end{small}
\end{table}
%% ------------------------------------------------------------------------- %%

\subsection{ADI splitting}

The ADI method is a predictor-corrector scheme as a typical example of algebraic splitting, which happens after the  discretization of equations. Here we are considering the discretized ODE version of the Peaceman-Rachford scheme \cite{Peachman} for (\ref{eqn:split-ode}). When $\Lambda = 2$, the ADI scheme takes the form

\begin{equation}
   \label{3.3ADI}
\left\{
  \begin{array}{ll}
  \vspace{0.05in}
    \displaystyle{ \frac{\widetilde{\upsilon}-\upsilon^n}{\frac{1}{2}\Delta t}} = f_1(t_{n+\frac{1}{2}} ,\widetilde{\upsilon}) + f_2(t_n, {\upsilon}^n), & \hbox{} \\
     \displaystyle{ \frac{\upsilon^{n+1}-\widetilde{\upsilon}}{\frac{1}2{}\Delta t}}  =  f_1(t_{n+\frac{1}{2}} ,\widetilde{\upsilon}) + f_2(t_{n+1}, {\upsilon}^{n+1}) . & \hbox{}
  \end{array}
\right.
\end{equation}

%% ------------------------------------------------------------------------- %%
The Butcher tableau for the scheme (\ref{3.3ADI}) is shown in Table \ref{table6},
and clearly, we see that we can view the ADI splitting scheme as a 2-stage
ARK method.
%% ------------------------------------------------------------------------- %%
%
%
%% ------------------------------------------------------------------------- %%
\begin{table}[!htb]
\centering
\bigskip
%\begin{small}
\begin{tabular}{c | ccc | cc c  }
 0 &  0      &   0    &     0     &    0  & 0 & 0     \\[4pt]
 $\frac{1}{2}$ &  0      &   $\frac{1}{2}$    &     0     &    $\frac{1}{2}$  & 0 & 0     \\[4pt]
 1 &  0      &   1    &     0     &    $\frac{1}{2}$   & 0 & $\frac{1}{2}$      \\[4pt]
 \hline
   &  0      &   1    &     0     &    $\frac{1}{2}$   & 0 &   $\frac{1}{2}$  \\
 %\hline
   \end{tabular}
   \caption{Butcher tableau for ADI scheme.}
   \label{table6}
%\end{small}
\end{table}
%% ------------------------------------------------------------------------- %%

% section 3
%\input{section3.tex}

\section{Formulation of IDC-OS schemes}

In this section, we review the formulation of IDC-OS
%the IDC-OS scheme
presented in \cite{Andrew2}. The authors \cite{Andrew2} considered
IDC methods for implicit-explicit (IMEX) schemes, where the non stiff part of the problem was
treated explicitly, and the stiff part of the problem was treated implicitly.  At
present, our focus is on entirely implicit schemes.

We begin with some preliminary definitions. The starting point is to
partition the time interval $[0, T]$ into intervals
$[t_n, t_{n+1}]$, $ n=0, 1, ..., N-1$, that satisfy
\begin{eqnarray}
\label{eqn_idc_os_2}
    0 = t_0 < t_1 < t_2 < \cdots < t_n < \cdots < t_N = T.
\end{eqnarray}
``macro''-time steps are defined by $H_n = t_{n+1} - t_n$, and we permit them to
vary with $n$.
Next, each interval $[t_n, t_{n+1}]$ is further partitioned into M sub-intervals
$[t_{n,m}, t_{n,m+1}]$, $m = 0,1,..., M-1$,
\begin{eqnarray}
\label{eqn_idc_os_3}
t_{n} = t_{n, 0} < t_{n, 1} < t_{n, 2} < \cdots < t_{n, m} < \cdots < t_{n, M} = t_{n+1}
\end{eqnarray}
with time step size $h_{n, m} = t_{n, m} - t_{n, m-1}$. If Gaussian quadrature
nodes are selected, as was originally done with the SDC method
\cite{Dutt}, $h_{n, m}$ varies with $m$. Here, we only consider the case of
uniform quadrature nodes, i.e. with $h_{n, m } = \frac{H_n}{M}$ for $m = 1, 2,
\dots, M$. Thus, without any ambiguity, we will drop the subscript $m$ on $h_{n, m}$. Note that in what follows we will use superscript $[i]$ to denote the $i^{th}$ correction at a discrete set of time points and superscript $(i)$ to denote the continuous approximation given by passing a $M^{th}$ order polynomial  through the discrete approximation. For simplicity, we drop the $n$ subscript for the description of the IDC procedure on ``macro''-time interval $[t_n, t_{n+1}]$. The whole iterative prediction-correction procedure is completed before moving on to the next time interval $[t_{n+1}, t_{n+2}]$. The numerical solution at $t_{n+1}$ serves as the initial condition for the following interval $[t_{n+1}, t_{n+2}]$.
\begin{itemize}
  \item \textbf{Prediction step :}
   Use an $r_0$-th order numerical method to obtain a preliminary solution to IVP (\ref{eqn:split-ode})
   \begin{eqnarray}\label{eqn_idc_os_4}
   \upsilon^{[0]} = (\upsilon^{[0]}_0, \upsilon^{[0]}_1, \dots, \upsilon^{[0]}_m, \dots, \upsilon^{[0]}_M),
   \end{eqnarray}
  which is an $r_0$-th order approximation to the exact solution
  \begin{eqnarray}\label{eqn_idc_os_5}
   u = (u_0, u_1, ..., u_m, ..., u_M),
   \end{eqnarray}
  where $u_m = u(t_m)$ is the exact solution at $t_m$ for $m = 0, 1, 2, ... , M$.

  \item \textbf{Correction step : }
  Use the error function to improve the accuracy of the scheme at each
  iteration. For $k = 1$ to $c_{s}$, ($c_s$ is the number of correction steps):
  \begin{description}
    \item(1) Denote the error function from the previous step as
    \begin{eqnarray}\label{eqn_idc_os_7}
    e^{(k-1)} (t) = u(t) - \upsilon^{(k-1)}(t),
    \end{eqnarray}
    where $u(t)$ is the exact solution and $\upsilon^{(k-1)}(t)$ is an $M$-th degree polynomial interpolating $\upsilon^{[k-1]}$. Note that the error function, $e^{(k-1)}(t)$ is not a polynomial in general.

    \item(2) Denote the residual function as
    \begin{eqnarray}\label{eqn_idc_os_8}
    \epsilon^{(k-1)}(t) \equiv (\upsilon^{(k-1)})'(t)- f(t, \upsilon^{(k-1)}) , %- L_2(t, \upsilon^{(k-1)})
    \end{eqnarray}
    and compute the integral of the residual. For example,
    \begin{eqnarray}\label{eqn_idc_os_9}
    \int^{t_{m+1}}_{t_0}\epsilon^{(k-1)}(\tau)  d \tau \approx \upsilon^{[k-1]}_{m+1} - u_0 -(t_{m+1}-t_0) \sum^M_{j=0}{\gamma}_{m, j}
     f(t_j, \upsilon^{[k-1]}_j),
 %   (L_1(t_j,\upsilon^{[k-1]}_j) + L_2(t_j, \upsilon^{[k-1]}_j))
    \end{eqnarray}
 { where ${\gamma}_{m, j}$ are the coefficients that result from approximation of the integral by quadrature formulas} and $\upsilon^{[k-1]}_j = \upsilon^{(k-1)}(t_j)$.

    \item(3) Use an $r_k$-th order numerical method to obtain an approximation to error vector
\begin{eqnarray}\label{eqn_idc_os_11}
    e^{[k-1]} = (e^{[k-1]}_0, ... , e^{[k-1]}_m, ... , e^{[k-1]}_M),
    \end{eqnarray}
    where $e^{[k-1]}_m = e^{(k-1)}(t_m)$ is the value of the exact error function (\ref{eqn_idc_os_7}) at time $t_m$
and we denote it as
    \begin{eqnarray}\label{eqn_idc_os_10}
    \delta^{[k]} = (\delta^{[k]}_0, ... , \delta^{[k]}_m, ... , \delta^{[k]}_M).
    \end{eqnarray}
    To compute $\delta^{[k]}$ by an operator splitting method consistent with the base method, we first express the error equation in a form consistent with original problem we are solving. We start by  differentiating the error (\ref{eqn_idc_os_7}), together with (\ref{eqn:split-ode})
\begin{eqnarray}\label{eqn_idc_os_11.5}
    (e^{(k-1)})'(t) & = & u'(t) - (\upsilon^{(k-1)})'(t) \\ \nonumber
                    & = & f(t, u(t)) - f(t, \upsilon^{(k-1)}(t)) - \epsilon^{(k-1)}(t) \\ \nonumber
                    & = & f(t, \upsilon^{(k-1)}(t)+e^{(k-1)}(t)) - f(t, \upsilon^{(k-1)}(t)) - \epsilon^{(k-1)}(t) . \nonumber
    \end{eqnarray}

    Bring the residual to the left hand side, we have
    \begin{equation}\label{eqn_idc_os_13}
    (e^{(k-1)}(t) + \int^t_{t_0}\epsilon^{(k-1)} (\tau) d\tau  )' =  f(t, \upsilon^{(k-1)}(t)+e^{(k-1)}(t)) - f(t,\upsilon^{(k-1)}(t)) .
    \end{equation}
   We now make the following change of variable,
    \begin{eqnarray}\label{eqn_idc_os_14}
   & Q^{(k-1)} (t) = e^{(k-1)}(t) + \displaystyle \int^t_{t_0} \epsilon^{(k-1)} (\tau) d\tau, \\ \nonumber
   & G^{(k-1)}(t, Q^{(k-1)}(t)) = f(t, \upsilon^{(k-1)}(t)+Q^{(k-1)}(t)-\displaystyle \int^t_{t_0} \epsilon^{(k-1)} (\tau) d\tau) - f(t,\upsilon^{(k-1)}(t)).
    \end{eqnarray}
    With this change of variable, we see that the  error equation can be
    expressed as an IVP of the form,
    \begin{align}
 \begin{cases}\label{eqn_idc_os_15}
 \vspace{0.1in}
 \displaystyle{(Q^{(k-1)})'(t)}=   G^{(k-1)}(t, Q^{(k-1)}(t)) , \qquad t\in [t_0,t_{M}],\\
                  %               \vspace{0.1in}
  Q^{(k-1)}(t_0) = 0.
  \end{cases}
\end{align}
This is now in the form of \eqref{eqn:split-ode} and we can apply the
same operator splitting scheme to (\ref{eqn_idc_os_15}) that we
applied to (\ref{eqn:split-ode}) and obtain the numerical approximation to ${\vartheta}^{[k-1]}_m = Q^{(k-1)}(t_m)$.  Recovering $\delta$ given ${\vartheta}$ is a simple procedure. \\

    \item(4) Update the numerical solution as $\upsilon^{[k]} = \upsilon^{[k-1]} + \delta^{[k]}$.

  \end{description}

\noindent {\bf Remark 1 (The prediction step):} For example, if we apply the discrete form of
first order Lie-Trotter splitting (\ref{3.3}) to (\ref{eqn:split-ode}) with $\Lambda = 2$, we have for $m = 0, 1, 2, ... , M-1$,

  \begin{equation}
   \label{eqn_idc_os_6}
\left\{
  \begin{array}{ll}
  \vspace{0.05in}
    \displaystyle{ \frac{\widetilde{\upsilon}-\upsilon^{[0]}_m}{h_n}} = f_1(t_{m+1},\widetilde{\upsilon}), & \hbox{} \\
     \displaystyle{ \frac{\upsilon^{[0]}_{m+1}-\widetilde{\upsilon}}{h_n}}  = f_2(t_{m+1}, \upsilon^{[0]}_{m+1}) . & \hbox{}
  \end{array}
\right.
\end{equation}

\noindent {\bf Remark 2 (The correction step):}
As an example, if we use ADI splitting in the correction step, we will  solve  (\ref{eqn_idc_os_15}) with $\Lambda = 2$, we have for $m = 0, 1, 2, ..., M-1$,
\begin{equation}
   \label{eqn_idc_os_16}
\left\{
  \begin{array}{ll}
  \vspace{0.05in}
    \displaystyle{\frac{\widetilde{\vartheta} - \vartheta^{[k]}_m}{\frac{h_n}{2}}}  =  G_1^{(k-1)}(t_m+\frac{h_n}{2}, \widetilde{\vartheta})  + G_2^{(k-1)}(t_m, \vartheta^{[k]}_m)  , & \hbox{} \\
     \displaystyle{\frac{\vartheta^{[k]}_{m+1} - \widetilde{\vartheta}}{\frac{h_n}{2}}  = G_2^{(k-1)}(t_{m+1}, \vartheta^{[k]}_{m+1}  ) + G_1^{(k-1)}(t_m+\frac{h_n}{2}, \widetilde{\vartheta} )}, & \hbox{}
  \end{array}
\right.
\end{equation}
where
\begin{equation}\label{eqn_idc_os_17}
G_{\nu}^{(k-1)}(t, Q^{(k-1)}(t)) = f_{\nu}(t, \upsilon^{(k-1)}(t)+Q^{(k-1)}(t)-\displaystyle \int^t_{t_0} \epsilon^{(k-1)} (\tau) d\tau) - f_{\nu}(t,\upsilon^{(k-1)}(t))
\end{equation}
for $\nu = 1, 2$. Moreover, we note that  we split the residual term equally for each operator in implementation.

\end{itemize}

% section 4
%\input{section4.tex}

\section{Analysis of IDC-OS methods}
\label{sec:analysis}

%\subsection{Truncation Error Analysis}
In this section, we will discuss the error estimate for IDC-OS
methods. Our analysis is similar to previous work of
IDC-RK and  IDC-ARK
\cite{Andrew, Andrew1, Andrew2}.
%% ------------------------------------------------------------------------- %%

%% ------------------------------------------------------------------------- %%
In section \ref{subsec:continuous-splitting-theory},
we will establish that the IDC procedure can successfully reduce the
splitting error for differential operator splitting methods where each
sub-problem is solved exactly. In section \ref{subsec:discrete-splitting-theory}, we continue by
leveraging the ideas from the work in \cite{Andrew2}, and prove that the
overall  accuracy for the fully discrete methods is increased, as expected, with
each successive correction.  The second set of arguments apply to the
discrete form of the differential operator splitting methods as well as the
algebraic operator splitting methods.
We present results for the stability regions of IDC-OS schemes in section \ref{subsec:stability}.
We remark that throughout this section, superscripts with a curly bracket $\{k\}$ denote the analytical functions related to solutions through differential splitting methods.
%% ------------------------------------------------------------------------- %%

\subsection{Splitting error: exact solutions to sub-problems}  %IDC-OS with Analytical Solutions for Sub-Problems}
\label{subsec:continuous-splitting-theory}

Differential operator splitting introduces a splitting error.  If each sub-problem
is solved exactly, the overall method only contains splitting error.  Our
starting point is to prove that IDC framework can reduce this splitting error. The primary result from this subsection is given by the following theorem.
\begin{thm}\label{thm0}
Assume $u(t)$ is the exact solution to IVP \eqref{eqn:split-ode}. Consider one time interval of an IDC method with $t \in [0, h]$. Suppose
Lie-Trotter splitting (\ref{3.2}) is used in the prediction step and the successive
$c_s$ correction steps,  and the sub-problems in each step are solved exactly. If $u(t)$ and $f_{\nu}$ are at least
$(c_s +3)$ differentiable, then the
splitting error is of order $\BigOh(h^{{c_s}+2})$ after $c_s$ correction steps.
\end{thm}
The proof of Theorem \ref{thm0} follows by induction from the following two lemmas: Lemma \ref{thm0_lemma0} for the prediction step and Lemma \ref{thm0_lemma1} for the correction steps respectively.
%% ------------------------------------------------------------------------- %%

%% ------------------------------------------------------------------------- %%
\begin{lem}\label{thm0_lemma0}
(Prediction step) Consider IVP \eqref{eqn:split-ode} on the
interval $t\in[0, h]$. If $u(t)$ and $f_{\nu}$ satisfy the
smoothness requirements
in Theorem \ref{thm0}, and $u^{\{0\}}(t)$ is the solution obtained by applying Lie-Trotter splitting (\ref{3.2}) to \eqref{eqn:split-ode}, and the followed sub-problems are solved exactly,
then the splitting error scales as
\begin{align*}
\| e^{(0)} \| = \| u(h)- u^{\{0\}}(h) \| \sim \BigOh(h^{2}), \qquad t \in [0, h].
\end{align*}
\end{lem}
The conclusion of Lemma \ref{thm0_lemma0} is simply a restatement of what the
local error of splitting methods measures.
The splitting error is $\BigOh(h^2)$ for first order
Lie-Trotter splitting \cite{McLachlan1}.
%% ------------------------------------------------------------------------- %%

%% ------------------------------------------------------------------------- %%
\begin{lem}\label{thm0_lemma1}
(Correction step) Assume $u(t)$ is the solution to IVP (\ref{eqn:split-ode}) on the interval $t\in[0, h]$.  Let $u(t)$, and $f_{\nu}$ satisfy the smoothness requirements in Theorem \ref{thm0}. For $k\leq c_s$, let $u^{\{k\}} (t)$ be  the solution after the prediction step and $k$-th correction step via Lie-Trotter splitting method in Theorem \ref{thm0}.
If $\| e^{(k-1)}\| \sim {\BigOh}(h^{k+1})$, then { $\| e^{(k)}\| \sim \BigOh(h^{k+2})$ after k correction steps}.
\end{lem}
%% ------------------------------------------------------------------------- %%

%% ------------------------------------------------------------------------- %%
\emph{Proof:} We show the proof with the simple case $\Lambda = 2$. We have the error equation (\ref{eqn_idc_os_15}) after prediction and $(k-1)$ correction steps. Use the Lie-Trotter splitting method (\ref{3.2}) to solve (\ref{eqn_idc_os_15}), we have
\begin{align}
 \begin{cases}\label{lemma1_eqn6}
 \vspace{0.1in}
 \displaystyle{(Q_1^{\{k-1\}}(t))'}=  G_1^{(k-1)}(t, Q_1^{\{k-1\}}(t)),  \qquad t\in [0,h],\\
                  %               \vspace{0.1in}
  Q_1^{\{k-1\}}(0) = 0,
  \end{cases}
  \end{align}
and
\begin{align}
 \begin{cases}\label{lemma1_eqn7}
 \vspace{0.1in}
 \displaystyle{(Q_2^{\{k-1\}}(t))'}=  G_2^{(k-1)}(t, Q_2^{\{k-1\}}(t)),  \qquad t\in [0,h],\\
                  %               \vspace{0.1in}
  Q_2^{\{k-1\}}(0) = Q_1^{\{k-1\}}(h).
  \end{cases}
  \end{align}
%where
%\begin{eqnarray}\label{lemma1_eqn7.1}
%G_{\nu}^{(k-1)}(t, Q^{(k-1)}(t)) = L_{\nu}(t, Q^{(k-1)}(t)-\int_0^t \epsilon^{(k-1)}(\tau) d \tau ), \qquad \nu = 1, 2
%\end{eqnarray}
with $G_{\nu}^{(k-1)}(t, Q^{(k-1)}(t))$ defined in (\ref{eqn_idc_os_17}).
Hence $Q_2^{\{k-1\}}(h)$ is the approximation of $Q^{(k-1)}(h)$ solved by the Lie-Trotter splitting method. It's easy to see that
\begin{equation}\label{lemma1_eqn7.5}
e^{(k)}(h) = e^{(k-1)}(h) - e^{\{k-1\}}(h) = Q^{(k-1)}(h)- Q_2^{\{k-1\}}(h),
\end{equation}
for $t \in [0, h]$. To prove $ Q^{(k-1)}(h)- Q_2^{\{k-1\}}(h) \sim {\BigOh}(h^{k+2}) $,
we examine the scaled variant
\begin{equation}\label{lemma1_eqn8}
\bar{Q}^{(k-1)}(t) = \frac{1}{h^k} Q^{(k-1)}(t).
\end{equation}
With this new notation, IVP (\ref{eqn_idc_os_15}) can be equivalently written as
\begin{align}
 \begin{cases}\label{lemma1_eqn9}
 \vspace{0.1in}
 \displaystyle{(\bar{Q}^{(k-1)}(t))'}=  \bar{G}^{(k-1)}(t, \bar{Q}^{(k-1)}(t)),  \qquad t\in [0,h],\\
                  %               \vspace{0.1in}
  \bar{Q}^{(k-1)}(0) = 0.
  \end{cases}
  \end{align}
with
\begin{equation}\label{lemma1_eqn10}
\bar{G}^{(k-1)}(t, \bar{Q}^{(k-1)}(t)) = \frac{1}{h^k} G^{(k-1)}(t, h^k\bar{Q}^{(k-1)}(t))~.~
\end{equation}
Using the Lie-Trotter splitting method to solve IVP (\ref{lemma1_eqn9}) will give us
\begin{align}
 \begin{cases}\label{lemma1_eqn11}
 \vspace{0.1in}
 \displaystyle{(\bar{Q}_1^{\{k-1\}}(t))'}=  \bar{G}_1^{(k-1)}(t, \bar{Q}_1^{\{k-1\}}(t)),  \qquad t\in [0,h],\\
                  %               \vspace{0.1in}
  \bar{Q}_1^{\{k-1\}}(0) = 0,
  \end{cases}
  \end{align}
  and
\begin{align}
 \begin{cases}\label{lemma1_eqn12}
 \vspace{0.1in}
 \displaystyle{(\bar{Q}_2^{\{k-1\}}(t))'}=  \bar{G}_2^{(k-1)}(t, \bar{Q}_2^{\{k-1\}}(t)),  \qquad t\in [0,h],\\
                  %               \vspace{0.1in}
  \bar{Q}_2^{\{k-1\}}(0) = \bar{Q}_1^{\{k-1\}}(h)~,~
  \end{cases}
  \end{align}
  with
  \begin{equation}\label{lemma1_eqn12.1}
\bar{G}_{\nu}^{(k-1)}(t, \bar{Q}^{(k-1)}(t)) = \frac{1}{h^k} G_{\nu}^{(k-1)}(t, h^k\bar{Q}^{(k-1)}(t)), \qquad \nu = 1, 2.
\end{equation}
 $\bar{Q}_2^{\{k-1\}}(h)$ is the approximation to $\bar{Q}^{(k-1)}(h)$ through Lie-Trotter splitting. If $ e^{(k-1)} \sim \BigOh(h^{k+1})$, it is easy to verify that $ Q^{(k-1)}(t) \sim \BigOh(h^{k+1})$ and $G^{(k-1)}(t, Q^{(k-1)}(t)) \sim \BigOh(h^{k+1})$. Similar as the work of IDC-RK in \cite{Andrew}, one can further check that $\frac{d}{d t}\bar{Q}^{\{k-1\}}(t) \sim \BigOh(1)$ and $\bar{G}^{(k-1)}(t, \bar{Q}^{\{k-1\}}(t)) \sim \BigOh(1)$ . Therefore,
\begin{equation}\label{lemma1_eqn13}
\parallel \bar{Q}^{(k-1)}(h)- \bar{Q}_2^{\{k-1\}}(h) \parallel \sim \BigOh(h^2).
\end{equation}
{Notice that IVP (\ref{lemma1_eqn6}) and (\ref{lemma1_eqn11}) are both first order  ODEs, and $h^k\bar{G}_1^{(k-1)}(t, \bar{Q}_1^{\{k-1\}}(t))  = G_1^{(k-1)}(t, h^k\bar{Q}_1^{\{k-1\}}(t))$. Since $\bar{Q}_1^{\{k-1\}}(t)$ is the solution to (\ref{lemma1_eqn11}),  $h^k\bar{Q}_1^{\{k-1\}}(t)$ is a solution to (\ref{lemma1_eqn6}). Through the uniqueness of the solution for IVP, one can conclude that
\begin{equation}\label{lemma1_eqn13.5}
\bar{Q}_1^{\{k-1\}}(h) = \frac{1}{h^k} Q_1^{\{k-1\}}(h)  .
%\qquad   \qquad \bar{Q}_2^{(k-1)}(h) = \frac{1}{h^k} Q_2^{(k-1)}(h)~.~
\end{equation}
Similarly, from IVP (\ref{lemma1_eqn7}) and (\ref{lemma1_eqn12}), one can further conclude
\begin{equation}\label{lemma1_eqn13.5_1}
%\widetilde{Q}_1^{(k-1)}(h) = \frac{1}{h^k} Q_1^{(k-1)}(h)  .
 \bar{Q}_2^{\{k-1\}}(h) = \frac{1}{h^k} Q_2^{\{k-1\}}(h).
\end{equation}
}
Thus (\ref{lemma1_eqn13}) is equivalent to
\begin{align}
    \label{lemma1_eqn14}
    \parallel \frac{1}{h^k} Q^{(k-1)}(h)- \frac{1}{h^k} Q_2^{(k-1)}(h) \parallel \sim \BigOh(h^2),
\end{align}
i.e.
\begin{align}\label{lemma1_eqn15}
    \parallel e^{(k)} \parallel = \parallel  Q^{(k-1)}(h)-  Q_2^{\{k-1\}}(h) \parallel \sim \BigOh(h^{k+2}).
\end{align}
We now complete the proof of Lemma \ref{thm0_lemma1} for the case of Lie-Trotter splitting.
\hfill $\square$
%% ------------------------------------------------------------------------- %%

%% ------------------------------------------------------------------------- %%
The conclusion in Theorem \ref{thm0} also holds for Strang splitting method (\ref{3.6}) and the proof is essentially the same as Lie-Trotter splitting. We have now demonstrated that IDC can lift the order of accuracy when each
sub-problem is solved exactly, however, in practice, we usually do not have
access to analytical solutions for these sub-problems. We will consider the fully discrete scheme in the next section.
%We now turn to proving
%that we can lift the overall accuracy for a fully discrete problem.

\subsection{Local truncation error: discrete solutions to sub-problems}
\label{subsec:discrete-splitting-theory}

%% ------------------------------------------------------------------------- %%
A fully discrete solution introduces additional error beyond the splitting
error.  In this section, we turn to analyzing fully discrete IDC-OS schemes and begin with some preliminary definitions \cite{Andrew}.
%% ------------------------------------------------------------------------- %%

%% ------------------------------------------------------------------------- %%
%% Definition
%% ------------------------------------------------------------------------- %%
\begin{de}
(Discrete differentiation) Consider the discrete data set,
$(\vec{t},\vec{\psi}) = \{(t_0, \psi_0),..., (t_M, \psi_M) \}$,
with $\{t_m\}^M_{m=0}$ defined as uniform quadrature nodes in (\ref{eqn_idc_os_3}).
We denote $L^M$ as the $M$-th degree Lagrangian interpolant of $(t, \psi)$:
\begin{align}
\label{eqn4.1}
    L^M(t, \psi) = \sum^M_{m=0} c_m(t)\psi_m, \qquad c_m(t) = \prod_{n \neq m} \frac{t-t_n}{t_m-t_n}.
\end{align}
An $s$-th degree discrete differentiation is a linear mapping that maps $\vec{\psi}$ to $\overrightarrow{\hat{d}_s \psi}$, where
\begin{eqnarray}
\label{eqn4.2}
(\hat{d}_s \psi)_m = \frac{\partial^s}{\partial t^s} L^M(t, \psi)\mid_{t=t_m}.
\end{eqnarray}
This linear mapping can be represented by a matrix multiplication $\overrightarrow{\hat{d}_s \psi} = \hat{D}_s \cdot \vec{\psi}$, where $\hat{D}_s \in \Re^{(M+1)\times{(M+1)}}$ and $(\hat{D})_{mn} = \frac{\partial^s}{\partial t^s} c_n(t)\mid_{t=t_m}$, $m,n = 0,...,M.$

\end{de}
%We notice that given a distribution of quadrature nodes on interval $[0,1]$, the differentiation matrices, $\hat{D}^{[0,1]}_s, s = 1,...,M$ contains constant entries. If this distribution of quadrature nodes is rescaled from $[0,1]$ to $[0,H]$, then the corresponding differentiation matrices are $\hat{D}_1 = \frac{1}{H}\hat{D}^{[0,1]}_1$ and $\hat{D}_s = (\frac{1}{H})^s \hat{D}^{[0,1]}_s$.

%% ------------------------------------------------------------------------- %%
%% Definition
%% ------------------------------------------------------------------------- %%
\begin{de}
The $(\hat{S}, \infty)$ Sobolev norm of the discrete data set $(\vec{t}, \vec{\psi})$ is defined as
\begin{eqnarray}
\label{eqn4.3}
\|\vec{\psi}\|_{\hat{S}, \infty} \doteq \sum^{\hat{S}}_{s=0} \parallel\overrightarrow{\hat{d}_s\psi}\parallel_{\infty} = \sum^{\hat{S}}_{s=0} \parallel \hat{D}_s \cdot \vec{\psi}\parallel_{\infty},
\end{eqnarray}
where $\overrightarrow{\hat{d}_s \psi} = Id \cdot \hat{\psi}$ is the identity matrix operating on $\hat{\psi}$.
\end{de}

%% ------------------------------------------------------------------------- %%
%% Definition
%% ------------------------------------------------------------------------- %%
\begin{de}
(smoothness of a discrete data set) A discrete data set $(\vec{t},\vec{\psi}) = \{(t_0, \psi_0),..., (t_M, \psi_M) \}$ possesses $\hat{S}  (\hat{S} \leq M)$ degrees of smoothness if $\parallel \vec{\psi}\parallel_{\hat{S},\infty}$ is bounded as $h\rightarrow 0$, with $h$ defined as the step size in the sub-interval $(t_{m},t_{m+1})$ where $m = 0,1,\cdots,M-1$.
\end{de}
%% ------------------------------------------------------------------------- %%

%% ------------------------------------------------------------------------- %%
As discussed in section 2, all the listed operator splitting schemes
are a form of ARK methods. Therefore, we can use the framework of the IDC-ARK
schemes in  {\cite{Andrew2}} to enhance the order of the discretized scheme.
%The difference between the work in { \cite{Andrew2}} and the work in this
%paper is that IDC-OS methods in this paper are fully implicit.
Hence,  we shall
describe only what is needed for clarity when   extending the results of the
work in {\cite{Andrew2}} to the fully implicit case under consideration
here.  For further details, we refer the reader to { \cite{Andrew1,
Andrew2}}.  The theorems below apply to lifting the order of algebraic splitting as well as the discrete form of differential splitting.
The splitting error discussed in Theorem \ref{thm0} is directly related to the local truncation error.  We note that the results in the following theorem can be generalized to all IDC-OS schemes which can be written as a form of ARK method and the proof is quite similar.
%% ------------------------------------------------------------------------- %%

%% ------------------------------------------------------------------------- %%
%% Theorem
%% ------------------------------------------------------------------------- %%
\begin{thm} \label{thm1}
Let $u(t)$ be the solution to IVP  {
(\ref{eqn:split-ode})}. Assume $u(t)$, $f(t,u)$ and $f_{\nu}(t,u)$ are at least $\sigma$ differentiable with respect to each argument, where $\sigma \geq M+2$. Consider one time interval of an IDC
method with $t\in [0,H]$ and $M+1$ uniformly distributed quadrature points.
Suppose an $r_0$-th order ARK method \eqref{eqn:ark} is used in the prediction step and $(r_1, r_2, ... , r_{c_s})$-th order ARK methods are used in the successive $c_s$ correction steps. Let $s_k = \Sigma^k_{j=0}r_j$. If $s_{c_s}\leq M+1$,
then the local truncation error is of order ${\BigOh}(h^{s_{c_s}+1})$ after $c_s$ correction steps.
\end{thm}

The proof of Theorem \ref{thm1} follows by induction from the following lemmas for the prediction and correction steps.  For clarity, similar as \cite{Andrew2}, we will sketch a proof for Lie-Trotter splitting.

\begin{lem}\label{lem1}
(prediction step) Consider an $r_0$-th order ARK method for (\ref{eqn:split-ode}) on $[0,H]$, with (M+1) uniformly distributed quadrature points. $u(t)$ and $f_{\nu}$  satisfy the smoothness requirement in Theorem \ref{thm1} and let $\upsilon^{[0]} = (\upsilon^{[0]}_0, \upsilon^{[0]}_1, ... \upsilon^{[0]}_m, ..., \upsilon^{[0]}_M  )$ be the numerical solution. Then,
\begin{description}
  \item(1) The error vector $e^{[0]} = u - \upsilon^{[0]}$ satisfies $\|e^{[0]} \|_{\infty} \sim \BigOh(h^{r_0+1}) $.
  \item(2) The rescaled error vector $\displaystyle{\bar{e}^{[0]} = \frac{1}{h^{r_0}} e^{[0]}}$ has $\min(\sigma-r_0, M)$ degrees of smoothness in the discrete sense.
\end{description}
\end{lem}

\emph{Proof:}  (1) is obvious. We will prove (2) next. We drop the superscript $[0]$ as there is no ambiguity. Applying
the discrete form of the Lie-Trotter splitting (\ref{3.3}) to IVP
(\ref{eqn:split-ode}) with $\Lambda = 2$, we have
\begin{align}
 \begin{cases}\label{lem3_1}
 \vspace{0.1in}
 \displaystyle \frac{\widetilde{\upsilon} - \upsilon_m}{h} = f_1(t_{m+1}, \widetilde{\upsilon}),\\
                  %               \vspace{0.1in}
 \displaystyle \frac{\upsilon_{m+1} - \widetilde{\upsilon}}{h} = f_2(t_{m+1}, \upsilon_{m+1}),
  \end{cases}
  \end{align}
i.e.
\begin{equation}\label{lem3_2}
\upsilon_{m+1} = \upsilon_m + hf_1(t_{m+1}, \widetilde{\upsilon}) + hf_2(t_{m+1}, \upsilon_{m+1}).
\end{equation}
Performing Taylor expansion of $f_1(t_{m+1},\widetilde{\upsilon})$ at $t = t_m$, we get
\begin{equation}\label{lem3_3}
\upsilon_{m+1} = \upsilon_m + hf_1(t_{m}, \upsilon_m)+ hf_2(t_{m+1}, \upsilon_{m+1}) + \sum^{\sigma-2}_{i=1} \frac{h^{i+1}}{i !} \frac{d^i f_1}{d t^i} (t_m, \upsilon_m) + \BigOh(h^{\sigma}),
\end{equation}
on the other hand, the exact solution satisfies
\begin{eqnarray}\label{lem3_4}
u_{m+1} & = & u_m + \int^{t_{m+1}}_{t_m} f_1(\tau, u(\tau)) d\tau + \int^{t_{m+1}}_{t_m} f_2(\tau, u(\tau)) d\tau \\ \nonumber
    & = & u_m + h f_1(t_m, u_m) + \sum^{\sigma-2}_{i=1} \frac{h^{i+1}}{(i+1)!} \frac{d^i f_1}{d t^i}(t_m, u_m) \\ \nonumber
    & + &  h f_2(t_{m+1}, u_{m+1}) + \sum^{\sigma -2}_{i=1}  \frac{(-1)^{i+1} h^{i+1}}{(i+1)!} \frac{d^i f_2}{d t^i} (t_{m+1}, u_{m+1}) + \BigOh(h^{\sigma}).
\end{eqnarray}
Subtracting {(\ref{lem3_3}) }from (\ref{lem3_4}) gives
\begin{eqnarray*}\label{lem3_5}
e_{m+1} & = & e_m + h(f_1(t_m, u_m) - f_1(t_m, \upsilon_m)) + h (f_2(t_{m+1}, u_{m+1}) - f_2(t_{m+1}, \upsilon_{m+1})) \\ \nonumber
 & + & \sum^{\sigma-2}_{i=1} \frac{h^{i+1}}{(i+1)!} \frac{d^i f_1}{d t^i} (t_m, u_m) + \sum^{\sigma -2}_{i=1}  \frac{(-1)^{i+1} h^{i+1}}{(i+1)!} \frac{d^i f_2}{d t^i} (t_{m+1}, u_{m+1}) - \sum^{\sigma-2}_{i=1} \frac{h^{i+1}}{i !} \frac{d^i f_1}{d t^i} (t_m, \upsilon_m)   \\ &+& {\BigOh}(h^{\sigma}),
\end{eqnarray*}
where $e_{m+1} = u_{m+1}-\upsilon_{m+1}$ is the error at $t_{m+1}$. Denote
\begin{eqnarray}\label{lem3_6}
l_m = (f_1(t_m, u_m) - f_1(t_m, \upsilon_m)) +  (f_2(t_{m+1}, u_{m+1}) - f_2(t_{m+1}, \upsilon_{m+1}))
\end{eqnarray}
and
\begin{equation}\label{lem3_7}
r_m = \sum^{\sigma-2}_{i=1} \frac{h^{i+1}}{(i+1)!} \frac{d^i f_1}{d t^i} (t_m, u_m) + \sum^{\sigma -2}_{i=1}  \frac{(-1)^{i+1} h^{i+1}}{(i+1)!} \frac{d^i f_2}{d t^i} (t_{m+1}, u_{m+1}) - \sum^{\sigma-2}_{i=1} \frac{h^{i+1}}{i !} \frac{d^i f_1}{d t^i} (t_m, \upsilon_m).
\end{equation}
We will use an inductive approach with respect to the degree of the smoothness $s$ to investigate the smoothness of the rescaled error vector $\bar{e} = \frac{e}{h}$, and
\begin{equation}\label{lem3_8}
(d_1\bar{e})_m = \frac{\bar{e}_{m+1}-\bar{e}_m}{h} = \frac{l_m}{h} + \frac{r_m}{h^2} + {\BigOh}(h^{\sigma-2}).
\end{equation}
First of all, $\bar{e}$ has at least zero degrees of smoothness in the discrete sense since $\|\bar{e}\| \sim \BigOh(h)$. Assume $\bar{e}$ has $s\leq M-1$ degrees of smoothness, we will show $d_1\bar{e}$ has $s$ degrees of smoothness, from which we can conclude $\bar{e}$ has $(s+1)$ degrees of smoothness.
\begin{eqnarray}\label{lem3_9}
l_m &=& (f_1(t_m, u_m)-f_1(t_m, \upsilon_m)) + (f_2(t_{m+1}, u_{m-1}) - f_2(t_{m+1} \upsilon_{m+1})), \\ \nonumber
    &=& \sum^{\sigma-2}_{i=1} \frac{1}{i !} (e_m)^i \frac{\partial^i f_1}{\partial u^i} (t_m, u_m) + \sum^{\sigma-2}_{i=1} \frac{1}{i !} (e_{m+1})^i \frac{\partial^i f_2}{\partial u^i}(t_{m+1},u_{m+1}) \\ \nonumber & +& {\BigOh}((e_m)^{\sigma-1}) + {\BigOh}((e_{m+1})^{\sigma-1}) \\ \nonumber
    &=& \sum^{\sigma-2}_{i=1} \frac{h^i}{i !} (\bar{e}_m)^i \frac{\partial^i f_1}{\partial u^i} (t_m, u_m) + \sum^{\sigma-2}_{i=1} \frac{h^i}{i !} (\bar{e}_{m+1})^i \frac{\partial^i f_2}{\partial u^i}(t_{m+1},u_{m+1})\\ \nonumber & + & {\BigOh}((h \bar{e}_m)^{\sigma-1}) + {\BigOh}((h \bar{e}_{m+1})^{\sigma-1}).
\end{eqnarray}
By assuming that $f_1$ and $f_2$ have at least $\sigma$ degrees of smoothness, we can conclude $\frac{\partial^i f_1}{\partial u^i}$ and $\frac{\partial^i f_2}{\partial u^i}$ have at least $\sigma - i -1$ degrees of smoothness, which implies $h^{i-1}\frac{\partial^i f_1}{\partial u^i}$ and $h^{i-1}\frac{\partial^i f_2}{\partial u^i}$ have at least $\sigma-2$ degrees of smoothness. Therefore $\frac{l_m}{h}$ will have $\min(\sigma-2, s)$ degrees of smoothness.
Also, $\frac{r_m}{h^2}$ will have at least $s$ degrees of smoothness. Therefore, $d_1 \bar{e}$ has $s$ degrees of smoothness. Therefore, $\bar{e}$ has $(s+1)$ degrees of smoothness. Notice that $\sigma \ge M+2$, we complete the inductive approach and conclude $\bar{e}$ has $M$ degrees of smoothness. \hfill $\square$

Before investigating the correction step for IDC-OS schemes, we describe some details for the error equations first.  Notice that the error equation after (k-1) correction steps has the form of { (\ref{eqn_idc_os_13})}, with the notation $Q^{(k-1)}(t)$, we actually implement the problem (\ref{eqn_idc_os_15}) on time interval $[t_m, t_{m+1}]$ via discrete Lie-Trotter splitting as follows
\begin{align}
 \begin{cases}\label{lem4_1}
 \vspace{0.1in}
 \displaystyle \frac{\widetilde{\vartheta} - \vartheta^{[k]}_m}{h} = G_1^{(k-1)}(t_{m+1}, \widetilde{\vartheta}),\\
                  %               \vspace{0.1in}
 \displaystyle \frac{\vartheta^{[k]}_{m+1} - \widetilde{\vartheta}}{h} = G_2^{(k-1)}(t_{m+1}, \vartheta^{[k]}_{m+1}),
  \end{cases}
  \end{align}
through which $\vartheta^{[k]}_{m+1}$ is updated. Furthermore, we can update $\delta^{[k]}_{m+1}$ by (\ref{eqn_idc_os_14}) and (\ref{eqn_idc_os_9}).
Similarly, if we apply Lie-Trotter splitting to the scaled error equation (\ref{lemma1_eqn9}) over the time interval $[t_m, t_{m+1}]$, we have
\begin{align}
 \begin{cases}\label{lem4_3}
 \vspace{0.1in}
 \displaystyle \frac{\widetilde{\bar{\vartheta}} - \bar{\vartheta}^{[k]}_m}{h} = \bar{G}_1^{(k-1)}(t_{m+1}, \widetilde{\bar{\vartheta}}),\\
                  %               \vspace{0.1in}
 \displaystyle \frac{\bar{\vartheta}^{[k]}_{m+1} - \widetilde{\bar{\vartheta}}}{h} = \bar{G}_2^{(k-1)}(t_{m+1}, \bar{\vartheta}^{[k]}_{m+1}),
  \end{cases}
  \end{align}
from which we obtain $\bar{\vartheta}^{[k]}_{m+1}$ and further $\bar{\delta}^{[k]}_{m+1}$.

\begin{lem}\label{lem2}
(correction step) Let $u(t)$ and $L_{\nu}$ satisfy the smoothness requirements in Theorem \ref{thm1}. Suppose $e^{[k-1]} \sim {\BigOh}(h^{s_{k-1}+1})$
and $\displaystyle {\bar{e}^{[k-1]} = \frac{1}{h^{s_{k-1}}} e^{[k-1]}}$ has $(M+1-s_{k-1})$ degrees of smoothness in the discrete sense after the $(k-1)$-th correction step. Then, after the $k$-th correction step using an $r_k$-th order ARK method and $k \leq k_{c_s}$,
\begin{description}
  \item(1) {$ \| e^{[k]} \|_{\infty} \sim {\BigOh}(h^{s_k+1}) $}.
  \item(2) The rescaled error vector $\displaystyle{ \bar{e}^{[k]} = \frac{1}{h^{s_k}} e^{[k]} }$ has $M+1-s_k$ degrees of smoothness in the discrete sense.
 \end{description}
\end{lem}
\emph{Proof:} The proof of Lemma \ref{lem2} is similar as Lemma \ref{lem1}, but more tedious. Similar as in \cite{Andrew2}, we outline the proof here and present the difference between the proof of IDC-OS and IDC-RK, IDC-ARK in Proposition \ref{prop1}, we refer the reader to \cite{Andrew1} for details.
\begin{enumerate}
  \item Substract the numerical error vector from the integrated error equation
  \begin{equation}\label{lem4_eqn1}
  e^{[k]}_{m+1} = e^{[k-1]}_{m+1} - \delta^{[k]}_{m+1}
  \end{equation}
  and make necessary substitution and expansion via the rescaled equations.
  \item Bound the error $e^{[k]}$ by an inductive approach.
 \end{enumerate}

%%%%%%%%%%%%%%%%%%%%%%%%%%%%%%%%%%%%%%%%%%%%%%%%%%%%%%%%%%%%%%%%%%%%%%%%%%%

The following proposition is about the equivalence of the rescaled error vector and unscaled error vectors for Lie-Trotter splitting. We remark that the proof of this proposition shows the difference of the proof between IDC-OS and IDC-RK in \cite{Andrew1}, IDC-ARK in \cite{Andrew2}.
\begin{prop}\label{prop1}
Consider a single step of an IDC scheme constructed with the Lie-Trotter splitting scheme for the error equation, assume the exact solution $u(t)$, and $L_{\nu}$ satisfies the smoothness requirement in Theorem \ref{thm1}, then for a sufficiently smooth error function $e^{(k-1)}(t)$, the difference between the Taylor series for the exact error $e^{(k-1)}(t_{m+1})$ and the numerical error $\delta^{[k]}_{m+1}$ is ${\BigOh}(h^{k+2})$ after $k$ correction steps.
\end{prop}
\emph{Proof:} Notice that the left and right hand side terms of the rescaled error equation (\ref{lemma1_eqn9}) is ${\BigOh}(1)$, applying the discrete form of Lie-Trotter splitting scheme (\ref{3.3}) to (\ref{lemma1_eqn9}) will result in
\begin{equation}\label{prop1_eqn1}
\bar{Q}^{[k-1]}_{m+1} - \bar{\vartheta} ^{[k]}_{m+1} \sim {\BigOh}(h^2).
\end{equation}
Since
\begin{equation}\label{prop1_eqn2}
\bar{Q}^{[k-1]}_{m+1} = \frac{1}{h^k} Q^{[k-1]}_{m+1} = \frac{1}{h^k} (e^{[k-1]}_{m+1} - \int^{t_{m+1}}_{t_m} \varepsilon^{(k-1)}(\tau) d \tau).
\end{equation}
The proof is complete if the following argument holds.
\begin{equation}\label{prop1_eqn3}
h^k \bar{\delta}^{[k]}_{m} = \delta^{[k]}_{m} + {\BigOh}(h^{\sigma}), \qquad m = 0, 1, 2, ..., M,
\end{equation}
which is also equivalent to
\begin{equation}\label{prop1_eqn4}
h^k \bar{\vartheta}^{[k]}_{m} = \vartheta^{[k]}_{m} +{\BigOh}(h^{\sigma}), \qquad m = 0, 1, 2, ..., M.
\end{equation}
We will prove { (\ref{prop1_eqn4})} by induction. (\ref{prop1_eqn4}) holds for $m=0$ since the initial condition for the error equation is set as $0$. Assume (\ref{prop1_eqn3}) holds for $m$, then
\begin{eqnarray}\label{prop1_eqn5}
 \widetilde{\bar{\vartheta}} & = & \bar{\vartheta}^{[k]}_m + h \bar{G}^{(k-1)}_1 (t_{m+1}, \widetilde{\bar{\vartheta}} ), \\ \nonumber
 & = & \bar{\vartheta}^{[k]}_m + h  \sum^{\sigma -1}_{i = 0} \frac{h^i}{i !} \frac{d^i }{d t^i} \bar{G}_1^{(k-1)}(t_m, \bar{\vartheta}^{[k]}_m) + {\BigOh}(h^{\sigma}) \\ \nonumber
 & = & \bar{\vartheta}^{[k]}_m + h  \sum^{\sigma -1}_{i = 0} \frac{h^i}{i !} \frac{d^i }{d t^i} \left(\frac{1}{h^k} G_1^{(k-1)}(t_m, h^k \bar{\vartheta}^{[k]}_m) \right ) + {\BigOh}(h^{\sigma}) \\ \nonumber
 & = & \frac{1}{h^k} \left( \vartheta^{[k]}_m  +  h  \sum^{\sigma -1}_{i = 0} \frac{h^i}{i !} \frac{d^i }{d t^i}  G_1^{(k-1)}(t_m, h^k \bar{\vartheta}^{[k]}_m)  \right ) + {\BigOh}(h^{\sigma}).
\end{eqnarray}
On the other hand, Taylor expanding $\widetilde{\vartheta}$ at $t_m$ will give us
\begin{eqnarray}\label{prop1_eqn6}
\widetilde{\vartheta} & = & \vartheta^{[k]}_m + h G_1^{(k-1)}(t_{m+1}, \widetilde{\vartheta} )  \\ \nonumber
& = &  \vartheta^{[k]}_m  +  h  \sum^{\sigma -1}_{i = 0} \frac{h^i}{i !} \frac{d^i }{d t^i}  G_1^{(k-1)}(t_m, h^k \bar{\vartheta}^{[k]}_m) + {\BigOh}(h^{\sigma}).
\end{eqnarray}
Compare (\ref{prop1_eqn5}) and (\ref{prop1_eqn6}), we can conclude
\begin{equation}\label{prop1_eqn7}
\widetilde{\vartheta} = h^k \widetilde{\bar{\vartheta}} +{\BigOh}(h^{\sigma}).
\end{equation}
Similar approach to the second equation in (\ref{lem4_1}) and (\ref{lem4_3}) will result in
\begin{equation}\label{prop1_eqn8}
\vartheta^{[k]}_{m+1} =h^k\bar{\vartheta}^{[k]}_{m+1} + \BigOh(h^{\sigma}),
\end{equation}
which completes the inductive proof for (\ref{prop1_eqn4}). \hfill $\square$

%%%%%%%%%%%%%%%%%%%%%%%%%%%%%%%%%%%%%%%%%%%%%%%%%%%%%%%%%%%%%%%%%%%%%%%%%%%%%%%%%%

%% ------------------------------------------------------------------------- %%
%%
%% ------------------------------------------------------------------------- %%
\subsection{Stability}
\label{subsec:stability}
In this subsection, we study the linear stability of the proposed
IDC-OS numerical schemes.
As is common practice \cite{Hairer}, we consider the test problem
\begin{align}
\label{eqn6.1}
  u_t = \lambda u,
\end{align}
and observe how the numerical scheme behaves for different complex values of
$\lambda$.  Without loss of generality, we will assume that $u(0) = 1$, and
we'll consider a single time step of length $\Delta t = 1$.
The stability region of a numerical method is then defined as

\begin{align}
    \mathbb{D} := \{ \lambda \in \mathbb{C} : \left| u\left( 1 \right) \right| \leq 1 \}.
\end{align}
An additional complication comes from the fact that an operator splitting scheme
requires a splitting of the right hand side of \eqref{eqn6.1} into $\Lambda$
parts.  For simplicity,
we'll consider the special case of $\Lambda = 2$ with $\lambda = \lambda_1+\lambda_2$,
and we further assume that $\lambda_1 = \lambda_2$ for simplicity.
%% ------------------------------------------------------------------------- %%

\begin{figure}[!htb]
  %\centering
  \begin{minipage}[b]{0.3\textwidth}
    \centerline{
    \includegraphics[width=2.0in,angle=0,scale=1.0]{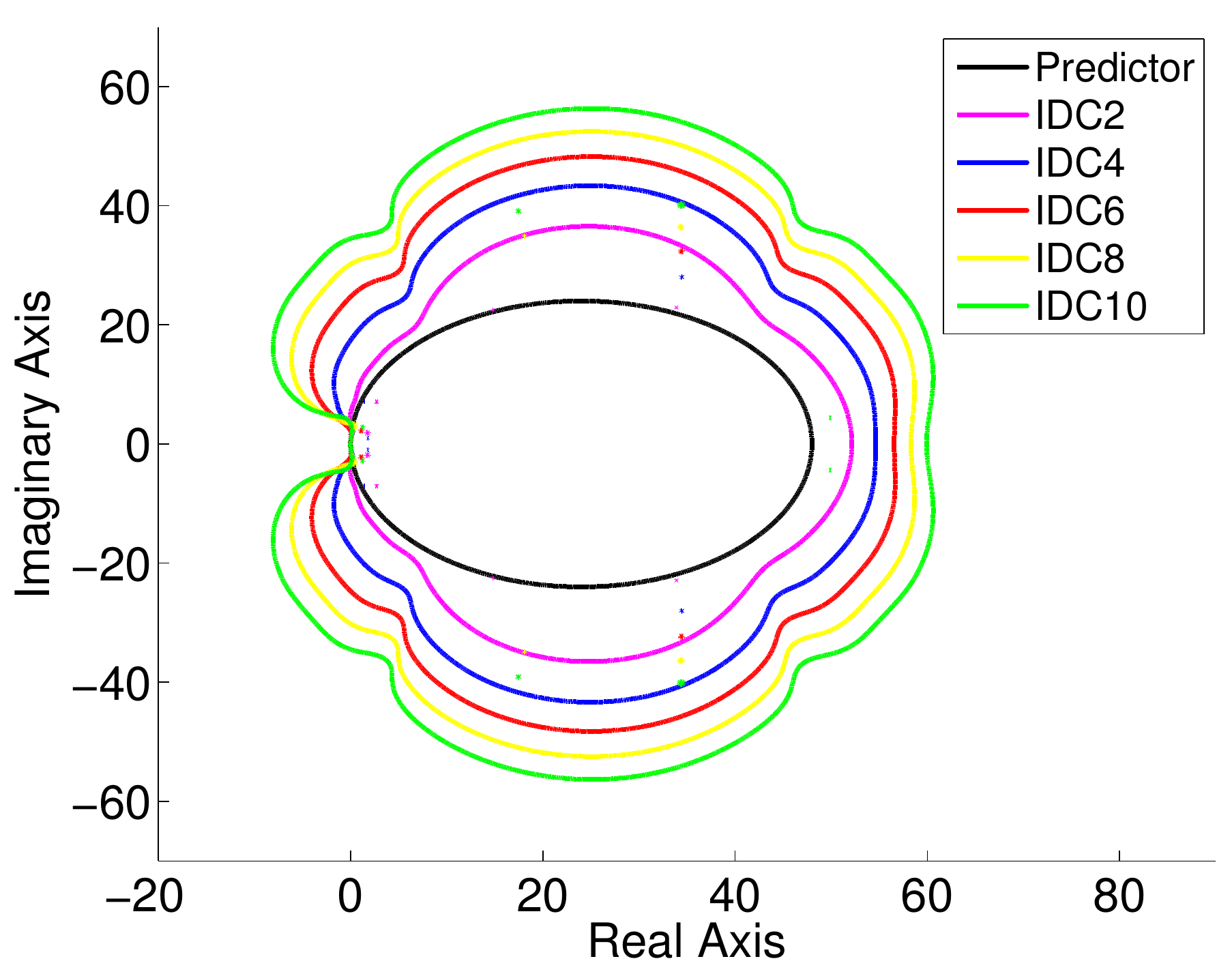}}
   % \caption{Caption 3}
     \qquad (a)
    \medskip
  \end{minipage}%
  \hspace{0.03\linewidth}%
  \begin{minipage}[b]{0.3\textwidth}
    \centerline{
    \includegraphics[width=2.0in,angle=0,scale=1.0]{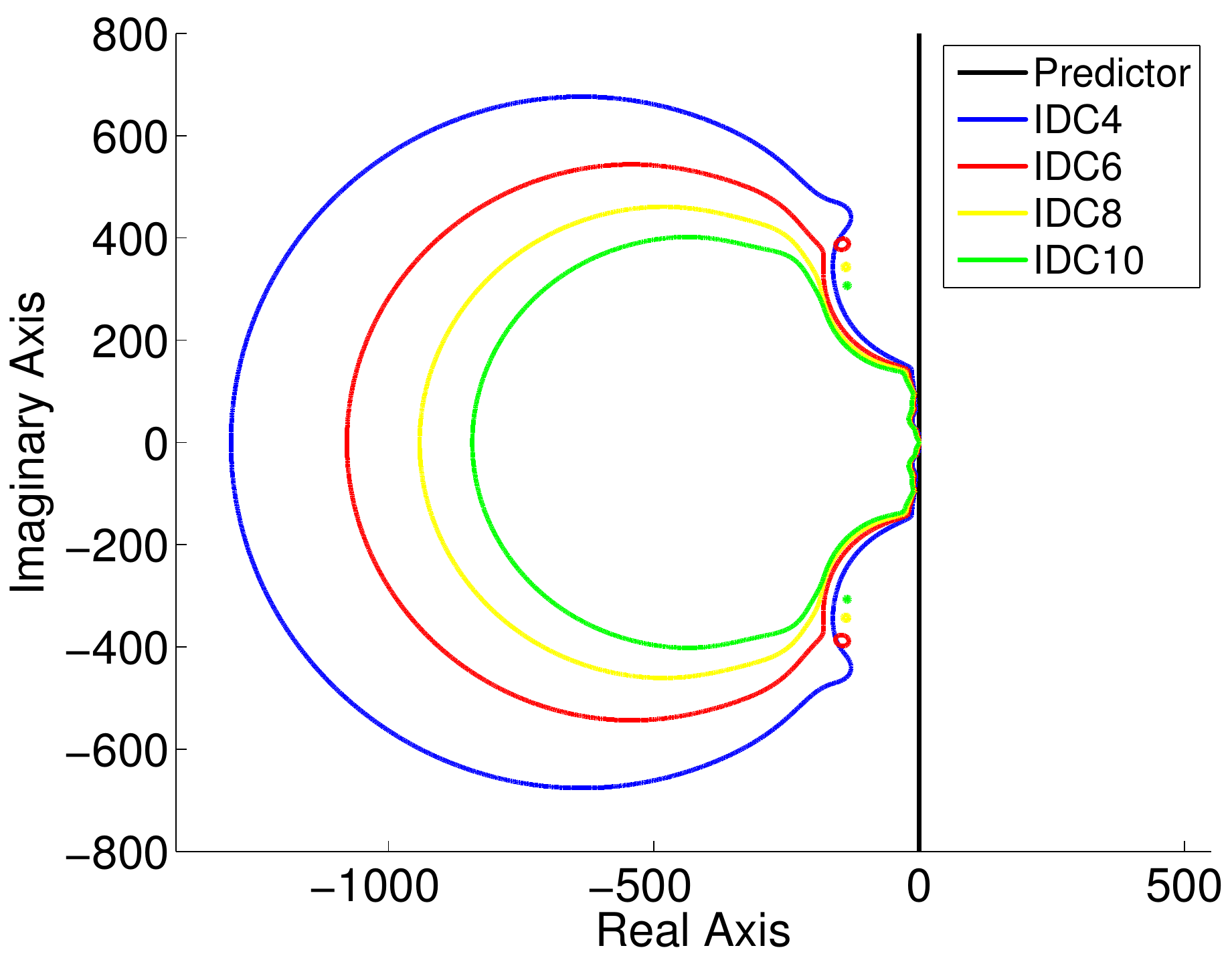}}
  %  \caption{Caption 1}
  \qquad (b)
  \medskip
  \end{minipage}
  \hspace{0.03\linewidth}%
  \begin{minipage}[b]{0.3\textwidth}
    \centerline{
    \includegraphics[width=2.0in,angle=0,scale=1.0]{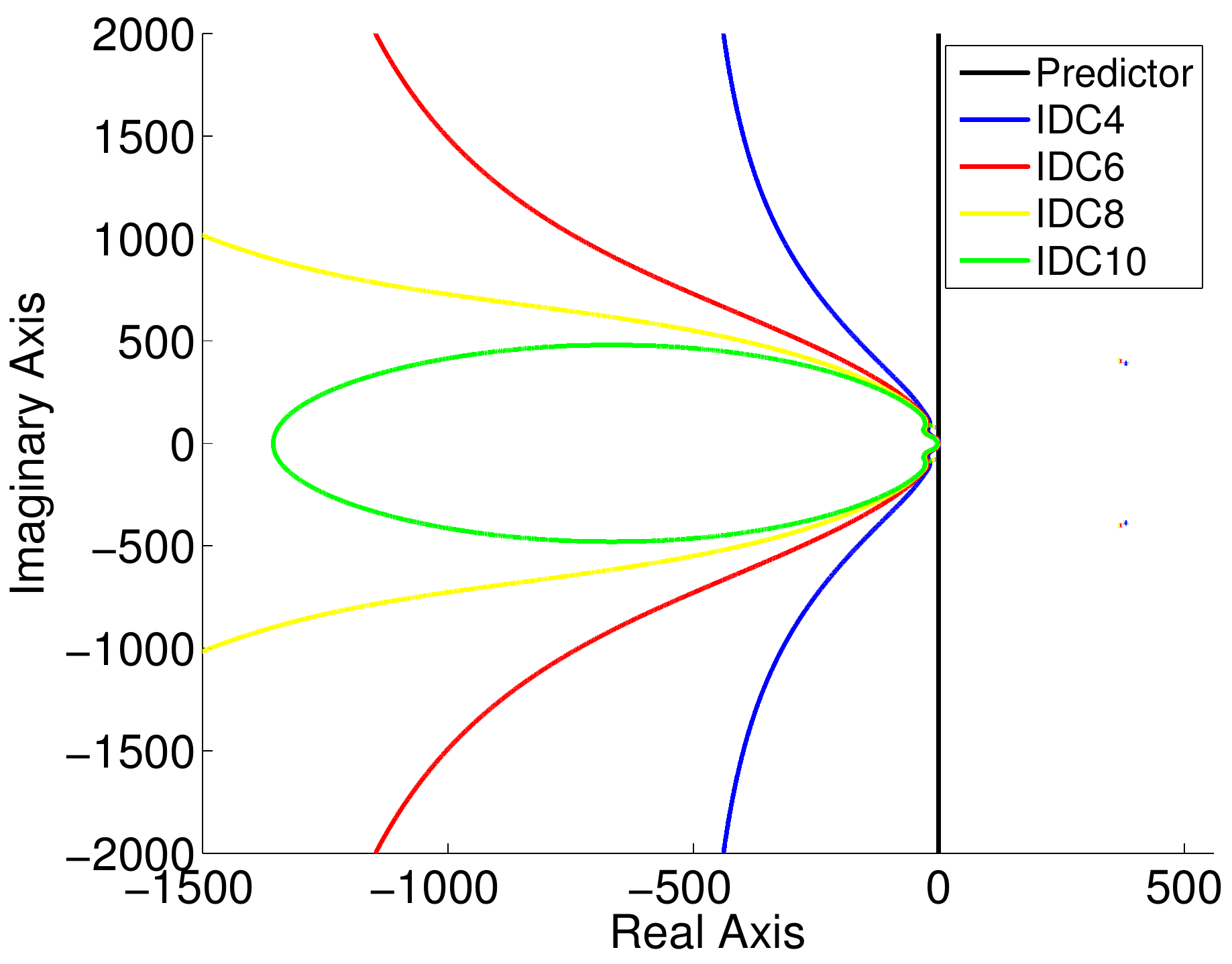}}
  %  \caption{Caption 1}
  \qquad (c)
  \medskip
  \end{minipage}
      \caption{Stability region for IDC-OS schemes with different number of corrections. (a) Lie-Trotter splitting; (b) Strang splitting; (c) ADI splitting.
 } \label{fig6_1}
\end{figure}

%% ------------------------------------------------------------------------- %%
In Figure \ref{fig6_1} we present stability regions for IDC-OS methods based on
three separate base solvers: Lie-Trotter splitting, Strang splitting and ADI
splitting. The stability region of Lie-Trotter splitting with IDC procedure is
everywhere outside the curves, and the stability regions for Strang splitting
and ADI is the finite region inside the curves.
The number in the legend denotes the order of the method. For example, ``IDC4" represents fourth order methods achieved by the IDC-OS schemes; for Lie-Trotter splitting, we require three correctors to attain fourth-order accuracy, whereas Strang and ADI splitting only require a single correction.
Our first observation is that of the three base solvers,
Lie-Trotter splitting is the only solver that retains an infinite region of
absolute stability, whereas Strang splitting and ADI reduce to finite regions of
absolute stability.

Consistent with other \emph{implicit} IDC methods,
the stability regions for our implicit
IDC-OS methods decreases as the number of correction
steps increases.
We have observed that larger stability regions can be found if we include more
quadrature nodes for evaluating the integral of the
residual.\footnote{For all simulations used in this work, we use 13 uniformly
distributed interior nodes for evaluating the residual integral in the error
equation.}
This leads us to conjecture
that a more accurate numerical approximation of the residual integral is
important in finding larger stability regions.

% section 5
%\input{section5.tex}
\section{Application of IDC-OS schemes to parabolic PDEs}
\label{sec:5}

In this section, we will discuss how to apply the IDC-OS framework to the parabolic problem of the form
\begin{align}
\begin{cases}\label{eqn:parabolic}
 \displaystyle{u_t}=  \nabla \cdot (a(x, y)\nabla u) +{s(t,u)}, \qquad  (x,y) \in \Omega \\
                  %               \vspace{0.1in}
 u(0,x,y)=u_0(x,y), \qquad \\
  u = g, \qquad (x,y) \in \partial \Omega.
\end{cases}
\end{align}
The methods can be generalized to a high dimensional setting, but in this work we restrict our attention to two dimensions. For differential splitting methods, it is quite straightforward to apply IDC-OS schemes if we solve (\ref{eqn:parabolic}) via method of lines. One can obtain semi-discrete ODE systems which have the same form as (\ref{eqn:split-ode}) after spatial discretization. It is natural to assume one operator, say $L_1$ is related to the terms in $x-$direction, while $L_2$ is related to the terms in $y-$direction.  As for algebraic splitting, one major difficulty for applying the IDC-OS framework to PDEs is how to handle the boundary and initial conditions for the error equation. In the following context, we will introduce one ADI formulation which can effectively deal with those issues. For simplicity, we only discuss the case when there is no nonlinear source in \eqref{eqn:parabolic}, i.e. $s(t, u) = 0$. 

%% ------------------------------------------------------------------------- %%

Classical ADI starts by applying second-order Crank-Nicholson time
discretization to the continuous PDE \eqref{eqn:parabolic}, this process produces
a semi-discrete scheme
\begin{equation}
\label{eqn3.4.2}
\frac{u^{n+1}-u^n}{\Delta t} = \frac{ a }{2}(u^{n+1}_{xx}+u^n_{xx})+\frac{
a_x}{2}(u^{n+1}_x+u^n_x)+\frac{a}{2}(u^{n+1}_{yy}+u^n_{yy})+\frac{a_y}{2}(u^{n+1}_y+u^n_y),
\end{equation}
where  $\Delta t = t^{n+1}-t^n$ is the time step, $a = a(x,y)$,  $a_x = a_x(x,y)$ and  $a_y= a_y(x,y)$.
On a two dimensional structured mesh, we choose to use the central difference
approximation (of orders $2,4$ or $6$) for approximating the spatial operators
$\frac{\partial}{\partial x^2}$, $\frac{\partial}{\partial x}$, $\frac{\partial}{\partial y^2}$ and $\frac{\partial}{\partial y}$,
and we denote them by $A_x, B_x, A_y, B_y$, respectively.
If the spatial discretization is performed on an $N_x\times N_y$ grid, there are $N_x\times N_y$ equations in the form of (\ref{eqn3.4.2}). We denote $\Upsilon$ as the unknowns in vector form, then we can write these $N_x\times N_y$ equations into matrix multiplication where
boundary conditions are also incorporated,
\begin{eqnarray}
\label{eqn3.4.3}
%\frac{\upsilon^{n+1}-U^n}{k} = \frac{1}{2}(A_x\textcolor{red}{\upsilon}^{n+1}+A_x\textcolor{red}{\upsilon}^n)+\frac{1}{2}(A_y\textcolor{red}{\upsilon}^{n+1}+A_yU^n)+
%\frac{1}{2}(g^{n+1}_{A_x}+g^n_{A_x}+g^{n+1}_{A_y}+g^n_{A_y})
\frac{{\Upsilon}^{n+1}-{\Upsilon}^n}{\Delta t}& =&  \nonumber
\frac{a}{2}(A_x{\Upsilon}^{n+1}+A_x{\Upsilon}^n)+\frac{a_x}{2}(B_x{\Upsilon}^{n+1}+B_x{\Upsilon}^n) \\ \nonumber
& + &\frac{a}{2}(A_y{\Upsilon}^{n+1}+A_y{\Upsilon}^n)+\frac{a_y}{2}(B_y{\Upsilon}^{n+1}+B_y{\Upsilon}^n)\\
& + &  \frac{a}{2}(g^{n+1}_{A_x}+g^n_{A_x}) +\frac{a_x}{2}(g^{n+1}_{B_x}+g^n_{B_x}) \\ \nonumber
& + & \frac{a}{2}(g^{n+1}_{A_y}+g^n_{A_y})+ \frac{a_y}{2}(g^{n+1}_{B_y}+g^n_{B_y})~,~
\end{eqnarray}
where $g_{A_x}$, $g_{B_x}$, $g_{A_y}$, and $g_{B_y}$ are the boundary {terms}. Notice that, different
from \cite{Doug3}, we enforce the boundary
conditions strictly in the scheme. It is easy to verify that the
method given in { (\ref{eqn3.4.3})} is second  order accurate in time. Specifically, if we use six order central difference for spatial derivatives such as in the numerical simulations, the local truncation error of { (\ref{eqn3.4.3})} is ${\BigOh}(\Delta t \Delta x^6+\Delta t^3)$.
Denoting
%$\textcolor{red}{J_1}=\frac{k}{2}A_x, L_2 = \frac{k}{2}A_y$,
\begin{eqnarray}
& &\nonumber {J_1} = \frac{\Delta t}{2}(aA_x+a_xB_x), \\
& & {J_2} = \frac{\Delta t}{2}(aA_y+a_yB_y), \\
& & \nonumber S = \frac{a}{2}(g^{n+1}_{A_x}+g^n_{A_x}) +\frac{a_x}{2}(g^{n+1}_{B_x}+g^n_{B_x}) +  \frac{a}{2}(g^{n+1}_{A_y}+g^n_{A_y})+ \frac{a_y}{2}(g^{n+1}_{B_y}+g^n_{B_y})~,~
\end{eqnarray}
{ (\ref{eqn3.4.3})} is  equivalent to
\begin{equation}
\label{eqn3.4.4}
%(I-\textcolor{red}{J_1-J_2})U^{n+1} = (I+{J_1+J_2})U^n+\frac{k}{2}(g^{n+1}_{A_x}+g^n_{A_x}+g^{n+1}_{A_y}+g^n_{A_y})
(I-{J_1}-{J_2}){\Upsilon}^{n+1} = (I+{J_1}+{J_2}){\Upsilon}^n+\Delta t S.
\end{equation}
To set up an ADI scheme, we follow \cite{Doug3} by adding one term ${J_1J_2}{\Upsilon}^{n+1}$ to both sides of { (\ref{eqn3.4.4})}, which results in
\begin{equation}
\label{eqn3.4.5}
(I-{J_1}-{J_2}+{J_1J_2}){\Upsilon}^{n+1} = (I+{J_1+J_2+J_1J_2}){\Upsilon}^n+{J_1J_2}({\Upsilon}^{n+1}-{\Upsilon}^n)+\Delta t S.
\end{equation}
Then it is straightforward to factor {(\ref{eqn3.4.5})} as
\begin{equation}
\label{eqn3.4.7}
(I-{J_1})(I-{J_2}){\Upsilon}^{n+1} = (I+{J_1})(I+{J_2}){\Upsilon}^n+{J_1J_2}({\Upsilon}^{n+1}-{\Upsilon}^n) +\Delta t S.
\end{equation}
Let us consider the second term on the right hand side of  { (\ref{eqn3.4.7})}.
Observe that
\begin{equation}
\label{eqn3.4.8}
\Upsilon^{n+1} = \Upsilon^n +{\BigOh}(\Delta t),
\end{equation}
and that $J_1$ and $J_2$ both carry a $\Delta t$ in them, we see that the term
${J_1J_2}({\Upsilon}^{n+1}-{\Upsilon}^n)\sim {\BigOh}(\Delta t^3)$. Hence, the second term on the right hand side of { (\ref{eqn3.4.7})} is the same order as the truncation error, thus can be dropped. Therefore, the scheme  reduces to
\begin{equation}
\label{eqn3.4.9}
(I-{J_1})(I-{J_2}){\Upsilon}^{n+1} = (I+{J_1})(I+{J_2}){\Upsilon}^n+\Delta t S.
\end{equation}
To solve (\ref{eqn3.4.9}), a two-step method was proposed in \cite{Doug2,Peachman},
\begin{align}
 \begin{cases}\label{eqn3.4.10}
 \vspace{0.1in}
 (I-{J_1})\tilde{{\Upsilon}}^{n+\frac{1}{2}} = (I+{J_2}){\Upsilon}^n + \frac{\Delta t}{2}S, \qquad \text{x-sweep},\\
                  %               \vspace{0.1in}
  (I-{J_2}){\Upsilon}^{n+1} = (1+{J_1})\tilde{{\Upsilon}}^{n+\frac{1}{2}} + \frac{\Delta t}{2}S,  \qquad \text{y-sweep}.
  \end{cases}
\end{align}
%% ------------------------------------------------------------------------- %%
However, to be symbolically consistent, symmetric and suited for IDC method, we choose to split the boundary values $S$ in the following way,
\begin{align}
 \begin{cases}\label{eqn3.4.11}
 \vspace{0.1in}
 (I-{J_1})\tilde{{\Upsilon}}^{n+\frac{1}{2}} = (I+{J_2}){\Upsilon}^n + S_1, \qquad \text{x-sweep},\\
                  %               \vspace{0.1in}
  (I-{J_2}){\Upsilon}^{n+1} = (1+{J_1})\tilde{{\Upsilon}}^{n+\frac{1}{2}} + S_2,  \qquad
  \text{y-sweep},
  \end{cases}
  \end{align}
with boundary terms defined as
 \begin{align}
 \begin{cases}\label{eqn3.4.12}
 \vspace{0.1in}
 S_1 = \frac{\Delta t}{2}(ag^{n+1}_{A_x}+a_xg^{n+1}_{B_x}+ag^n_{A_y}+a_yg^n_{B_y}),\\
                  %               \vspace{0.1in}
  S_2 = \frac{\Delta t}{2}(ag^n_{A_x}+a_xg^n_{B_x}+ag^{n+1}_{A_y}+a_yg^{n+1}_{B_y}).
  \end{cases}
  \end{align}
%% ------------------------------------------------------------------------- %%
%

%% ------------------------------------------------------------------------- %%
It should also be pointed out that the boundary values $S$ are associated with
boundary functions $g$ at time $t = t_n$ and $t = t_{n+1}$, instead of
$t_{n+\frac{1}{2}}$, therefore there is no error introduced from intermediate
values $\tilde{\Upsilon}$ on the boundary. This is important for setting the boundary
conditions when solving the error equation of IDC when we combine the ADI scheme
with the IDC methodology. Because the Dirichlet boundary conditions of  (\ref{eqn:parabolic}) are exact and accounted for in the formulation of the prediction, therefore, the boundary terms will not show up in the correction steps.
%% ------------------------------------------------------------------------- %%

% section 6
%\input{numerical.tex}

\section{Numerical examples}
\label{numer}

%% ------------------------------------------------------------------------- %%
In this section, we present numerical results for the proposed implicit IDC-OS
schemes on a variety of examples of the parabolic initial-boundary value problem
\eqref{eqn:parabolic}, where our aim is to demonstrate the efficiency of the
proposed time-stepping methods.  We begin with two
linear examples of \eqref{eqn:parabolic}, and then present an example
of the heat equation with a nonlinear forcing term.
Our final two examples come from mathematical biology: the Fitzhugh-Nagumo
reaction-diffusion model and the Schnakenberg model.
%% ------------------------------------------------------------------------- %%

%% ------------------------------------------------------------------------- %%
Our present work is in two-dimensions, and every result is performed on
a square domain with a cartesian grid.  We solve \eqref{eqn:parabolic} using $6$-th order central difference for the spatial discretization in order that the temporal error is dominant in the measured numerical error.
%% ------------------------------------------------------------------------- %%

%% ------------------------------------------------------------------------- %%
%% Example 1
%% ------------------------------------------------------------------------- %%
\noindent {\bf Example 1.} \textbf{Linear example: Dirichlet boundary conditions.}
We solve initial boundary value problem (\ref{eqn:parabolic}) with constant
coefficient $a(x,y) = 1$ in the domain $[-1,1]\times[-1,1]$.
Initial condition is taken as $u_0(x,y) = (1-y)e^x$ and time dependent boundary conditions are ${g(x,y,t)} = (1-y)e^{t+x}$. Therefore, (\ref{eqn:parabolic}) has the exact solution $u(x,y,t)=(1-y)e^{t+x}$.  $N_{x,y} = N_x = N_y$ represents the number of spatial grids in $x$- and $y$-direction. $N_t$ is the time steps used in the time interval $[0, T]$ where $T$ is end time.  $c_s$ is the number of correction steps. $u$ is the exact solution and ${\upsilon}$ as the numerical solution. We solve Example 1 by first order Lie-Trotter splitting (\ref{3.3}), second order Strang splitting (\ref{3.2.1}) and ADI splitting (\ref{eqn3.4.11}) and all the splitting performed via dimensional fashion, the numerical error are shown in Table \ref{table5.1}, Table \ref{table5.2} and Table \ref{table5.4} respectively. We can clearly conclude that the schemes achieve the designed order with IDC methodology, i.e. with one more correction step, the order of the scheme increases by 1 for Lie-Trotter splitting, while the order of the scheme increases by 2 for Strang splitting and ADI splitting.

\begin{table}[!htb]
\centering
\bigskip
\begin{tabular}{|c|c c|cc|c c|c c|}
\hline
$N_{x, y} = 45$ & \multicolumn{8}{c|}{Number of time steps $N_t$ }  \\
\hline
Correction & $N_t = 60 $ & order & $N_t = 80$ & order & $N_t = 100$ & order & $N_t = 120$ & order  \\
\hline
$c_s$ = 0 & 1.53e-5 & -- & 1.15e-5 & 0.99 & 9.24e-6  & 0.98& 7.70e-6 &    1.00 \\
\hline
$c_s$ = 1 & 1.90e-7 &-- & 1.09e-7 & 1.93 & 7.08e-8 & 1.93 & 4.94e-8 &    1.97\\
\hline
$c_s$ = 2 & 3.10e-9& -- & 1.47e-9 & 2.59 & 8.06e-10 & 2.69 & 4.87e-10 &   2.76\\
\hline
\end{tabular}
\caption{Linear example with Dirichlet boundary conditions. Errors $\parallel u-{\upsilon}_{N_t}\parallel_{\infty}$ for Lie-Trotter splitting method, $T = 0.025$.}
\label{table5.1}
\end{table}

\begin{table}[!htb]
\centering
\bigskip
\begin{tabular}{|c|c c|cc|c c|c c|}
\hline
$N_{x, y} = 45$ & \multicolumn{8}{c|}{Number of time steps $N_t$ }  \\
\hline
Correction & $N_t = 60 $ & order & $N_t = 80$ & order & $N_t = 100$ & order & $N_t = 120$ & order  \\
\hline
$c_s$ = 0 & 3.02e-5 & -- & 1.69e-5 & 2.02 & 1.08e-5  & 2.01& 7.55e-6 &    1.96 \\
\hline
$c_s$ = 1 & 7.15e-7 &-- & 2.45e-7 & 3.72 & 1.04e-7 & 3.84 & 5.20e-8 &    3.80\\
\hline
$c_s$ = 2 & 3.64e-10& -- & 8.06e-11 & 5.24 & 2.28e-11 & 5.66 & 7.16e-12 &   6.35\\
\hline
%CM = 3 & 6.39e-11& -- & 2.32e-11 & 3.5218 & 1.08e-11 & 3.4265& 6.16e-12 &  3.0796  \\
%\hline
\end{tabular}
\caption{Linear example with Dirichlet boundary conditions. Errors $\parallel u-{\upsilon}_{N_t}\parallel_{\infty}$ for Strang splitting method, $T = 0.025$.}
\label{table5.2}
\end{table}

\begin{table}[!htb]
\centering
\bigskip
\begin{tabular}{|c|c c|cc|c c|c c|}
\hline
$N_{x, y} = 150$ & \multicolumn{8}{c|}{Number of time steps $N_t$ }  \\
\hline
Correction & $N_t = 60 $ & order & $N_t = 80$ & order & $N_t = 100$ & order & $N_t = 120$ & order  \\
\hline
$c_s$ = 0 & 3.68e-5 & -- & 2.07e-5 & 2.00 & 1.32e-5  & 2.00 & 9.20e-6 &    2.00 \\
\hline
$c_s$ = 1 & 4.49e-7 &-- & 2.08e-7 & 2.67 & 1.03e-7 & 3.13 & 5.07e-8 &    3.91\\
\hline
$c_s$= 2 & 1.18e-7 & -- & 1.69e-8 & 6.74 & 4.74e-9 & 5.71 & 1.59e-9 &   5.98\\
\hline
%CM = 3 & 6.39e-11& -- & 2.32e-11 & 3.5218 & 1.08e-11 & 3.4265& 6.16e-12 &  3.0796  \\
%\hline
\end{tabular}
\caption{Linear example with Dirichlet boundary conditions. Errors $\parallel u-{\upsilon}_{N_t}\parallel_{\infty}$ for IDC-OS based on ADI splitting, $T = 0.025$.}
\label{table5.4}
\end{table}

%\pagebreak

%% ------------------------------------------------------------------------- %%
%% Example 2
%% ------------------------------------------------------------------------- %%
\noindent {\bf Example 2.} \textbf{Linear example: periodic boundary conditions.}
In this example, we will solve (\ref{eqn:parabolic}) when $a(x,y) = 2 + 0.5\sin(\pi(4x+y))$, and the initial condition $u_0(x,y) = \sin(2\pi(x+y))$.
%With this setup, we can assume periodic boundary condition here.
We compute errors using the difference between two successive refinements:
%The error in this example is obtained from the difference of the numerical
%solution between the coarse mesh and refined mesh, i.e.
\begin{equation}\label{eqnerror}
\text{error} = \parallel {\upsilon}_{N_t} - {\upsilon}_{\frac{N_t}{2}}\parallel_{\infty},
\end{equation}
where $N_t$ describes the number of time steps.

Again, we present results using three splitting options:
Lie-Trotter, Strang and ADI splitting.
Convergence studies are presented in Tables \ref{table5.5}, \ref{table5.6} and \ref{table5.7}, and we can
also observe that the schemes achieve the designed order.
Note that the numerical error for the two correctors when $N_t = 320$ is
not reliable in the cases of Strang and ADI splitting because of precision limitation.

\begin{table}[!htb]
\centering
\bigskip
\begin{tabular}{|c|c c|cc|c c|c c|}
\hline
$N_{x, y} = 45$ & \multicolumn{8}{c|}{Number of time steps $N_t$ }  \\
\hline
Correction & $N_t = 40 $ & order & $N_t = 80$ & order & $N_t = 160$ & order & $N_t = 320$ & order  \\
\hline
$c_s$ = 0 & 4.65e-3 & -- & 2.35e-3 & 0.98 & 1.18e-3  & 0.99 & 5.94e-4 &    0.99 \\
\hline
$c_s$ = 1 & 1.85e-4 &-- & 5.68e-5 & 1.70 & 1.63e-5 & 1.80 & 4.44e-6 &    1.88 \\
\hline
$c_s$ = 2 & 3.47e-6 & -- & 6.55e-7 & 2.41 & 1.19e-7 & 2.46 & 1.88e-8 &   2.66   \\
\hline
\end{tabular}
\caption{Linear example with periodic boundary conditions.
Errors $\parallel {\upsilon}_{N_t} - {\upsilon}_{\frac{N_t}{2}}\parallel_{\infty}$  for Lie-Trotter splitting method, $T = 0.025$.}
\label{table5.5}
\end{table}

\begin{table}[!htb]
\centering
\bigskip
\begin{tabular}{|c|c c|cc|c c|c c|}
\hline
$N_{x, y} = 45$ & \multicolumn{8}{c|}{Number of time steps $N_t$ }  \\
\hline
Correction & $N_t = 40 $ & order & $N_t = 80$ & order & $N_t = 160$ & order & $N_t = 320$ & order  \\
\hline
$c_s$ = 0 & 5.24e-5 & -- & 1.31e-5 & 2.00 & 3.29e-6  & 1.99 & 8.22e-7 &    2.00 \\
\hline
$c_s$ = 1 & 3.30e-9 &-- & 2.06e-10 & 4.00 & 1.29e-11 & 4.00 & 8.04e-13 &    4.00 \\
\hline
$c_s$ = 2 & 5.80e-12 & -- & 4.90e-14 & 6.89 & 7.77e-16 & 5.98 & 1.11e-16 &   2.81   \\
\hline
\end{tabular}
\caption{Linear example with periodic boundary conditions.
Errors $\parallel {\upsilon}_{N_t} - {\upsilon}_{\frac{N_t}{2}}\parallel_{\infty}$  for Strang splitting method, $T = 0.025$.}
\label{table5.6}
\end{table}

\begin{table}[!htb]
\centering
\bigskip
\begin{tabular}{|c|c c|cc|c c|c c|}
\hline
$N_{x, y} = 200$ & \multicolumn{8}{c|}{Number of time steps $N_t$ }  \\
\hline
Correction & $N_t = 40 $ & order & $N_t = 80$ & order & $N_t = 160$ & order & $N_t = 320$ & order  \\
\hline
$c_s$ = 0 & 7.77e-5 & -- & 1.94e-5 & 2.00 & 4.85e-6  & 2.00 & 1.21e-6 &    2.00 \\
\hline
$c_s$ = 1 & 1.93e-8 &-- & 1.20e-9  & 4.00 & 7.52e-11 &4.00 & 4.70e-12 &    4.00 \\
\hline
$c_s$ = 2 & 1.43e-11 & -- & 2.23e-13 & 6.01 & 3.56e-15 & 5.96 & 1.04e-16 &   5.10  \\
\hline
\end{tabular}
\caption{Linear example with periodic boundary conditions.
Errors $\parallel {\upsilon}_{N_t} - {\upsilon}_{\frac{N_t}{2}}\parallel_{\infty}$  for ADI splitting method, $T = 0.05$.}
\label{table5.7}
\end{table}

%\pagebreak

%% ------------------------------------------------------------------------- %%
%% Example 3
%% ------------------------------------------------------------------------- %%
\noindent {\bf Example 3.}  \textbf{Nonlinear equation with Dirichlet boundary conditions. }
We now test the proposed IDC-OS methods on a nonlinear example of \eqref{eqn:parabolic} with
a known exact solution.

\begin{align}
 \begin{cases}\label{eqn5.9}
 \displaystyle{u_t}=  u_{xx}+u_{yy} -u^2 + e^{-2t}\cos^2(\pi x)cos^2(\pi y) +(2\pi^2-1)e^{-t}\cos(\pi x)\cos(\pi y),  \\
 u(0,x,y)=\cos(\pi x)\cos(\pi y),
  \end{cases}
\end{align}
on the domain $(x,y) \in [-1,1]\times[-1,1]$.
The exact solution to this problem is
$u(x,y,t) = e^{-t} \cos(\pi x)\cos(\pi y)$.  Given that we have an exact solution,
all our numerical tests use exact boundary conditions from this solution.
%% ------------------------------------------------------------------------- %%

%% ------------------------------------------------------------------------- %%
An IDC-OS solver for \eqref{eqn5.9} requires a definition for how the splitting
will be performed.
Here, we choose to split the problem into three pieces:
$L_1$ and $L_2$ are the same as linear case,
while $L_3$ contains the remaining non-linear terms,
\begin{equation}\label{eqn5.13}
L_3(t, u) = -u^2 + e^{-2t}\cos^2(\pi x)cos^2(\pi y) +(2\pi^2-1)e^{-t}\cos(\pi x)\cos(\pi y).
\end{equation}
We use Newton-iteration to solve the discretized version of
$u_t = L_3(t,u)$.
%% ------------------------------------------------------------------------- %%

%% ------------------------------------------------------------------------- %%
In Tables \ref{table5.3.1} and \ref{table5.3.2} we present results from applying
the IDC-OS method with Lie-Trotter and Strang splitting as the base solvers.
In each case, we can see the successful increase of order after
each correction: one in the case of Lie-Trotter splitting, and two in the case
of Strang splitting.
%% ------------------------------------------------------------------------- %%

In this work, we do not present results for IDC-OS methods based on
ADI splitting for non-linear problems due to their computational complexity.
The high-order differential operator splitting methods discussed here are
much simpler than what would arise from using even low-order ADI splitting.
%% ------------------------------------------------------------------------- %%

%% ------------------------------------------------------------------------- %%
\begin{table}[!htb]
\centering
\bigskip
\begin{tabular}{|c|c c|cc|c c|c c|}
\hline
$N_{x, y} = 45$ & \multicolumn{8}{c|}{Number of time steps $N_t$ }  \\
\hline
Correction & $N_t = 60 $ & order & $N_t = 80$ & order & $N_t = 100$ & order & $N_t = 120$ & order  \\
\hline
$c_s$ = 0 & 6.88e-3 & -- & 5.16e-3 & 1.00 & 4.13e-3  & 1.00 & 3.44e-3 &    1.00 \\
\hline
$c_s$ = 1 & 7.31e-4 &-- & 4.33e-4 & 1.82 & 2.87e-4 & 1.84 & 2.03e-4 &    1.90\\
\hline
$c_s$ = 2 & 1.60e-5& -- & 6.95e-6 & 2.90 & 3.59e-6 & 2.96 & 2.06e-6 &   3.05\\
\hline
\end{tabular}
\caption{Nonlinear example with Dirichlet boundary conditions. Errors $\parallel u-\upsilon_{N_t}\parallel_{\infty}$ for IDC-OS based on Lie-Trotter splitting, $T = 0.025$.}
\label{table5.3.1}
\end{table}

\begin{table}[!htb]
\centering
\bigskip
\begin{tabular}{|c|c c|cc|c c|c c|}
\hline
$N_{x, y} = 100$ & \multicolumn{8}{c|}{Number of time steps $N_t$ }  \\
\hline
Correction & $N_t = 60 $ & order & $N_t = 80$ & order & $N_t = 100$ & order & $N_t = 120$ & order  \\
\hline
$c_s$ = 0 & 9.21e-5 & -- & 5.20e-5 & 1.99 & 3.34e-5  & 1.98 & 2.33e-5 &    1.99 \\
\hline
$c_s$ = 1 & 3.04e-6 &-- & 8.96e-7 & 4.25 & 3.38e-7 & 4.37 & 1.56e-7 &    4.24\\
\hline
$c_s$ = 2 & 1.87e-8& -- & 3.22e-9 & 6.11 & 8.33e-10 & 6.06 & 2.84e-10 &   5.90\\
\hline
\end{tabular}
\caption{Nonlinear example with Dirichlet boundary conditions. Errors $\parallel u-\upsilon_{N_t}\parallel_{\infty}$ for IDC-OS based on Strang splitting, $T = 0.01$. }
\label{table5.3.2}
\end{table}
%% ------------------------------------------------------------------------- %%

\begin{figure}[!htb]
  %\centering
  \begin{minipage}[b]{0.3\textwidth}
    \centerline{
    \includegraphics[width=2.0in,angle=0,scale=1.0]{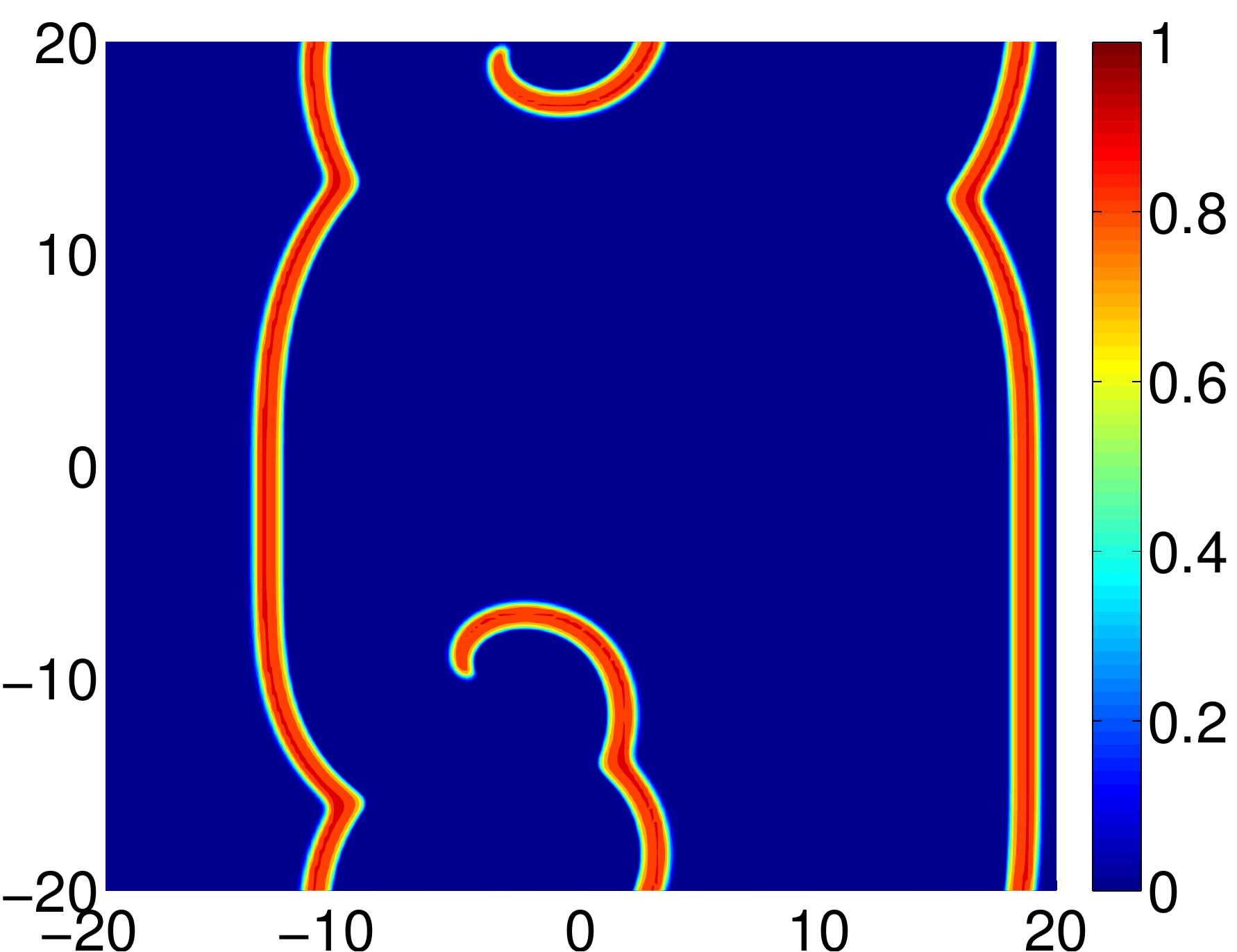}}
   % \caption{Caption 3}
      (a1)
    \medskip
  \end{minipage}%
  \hspace{0.03\linewidth}%
  \begin{minipage}[b]{0.3\textwidth}
    \centerline{
    \includegraphics[width=2.0in,angle=0,scale=1.0]{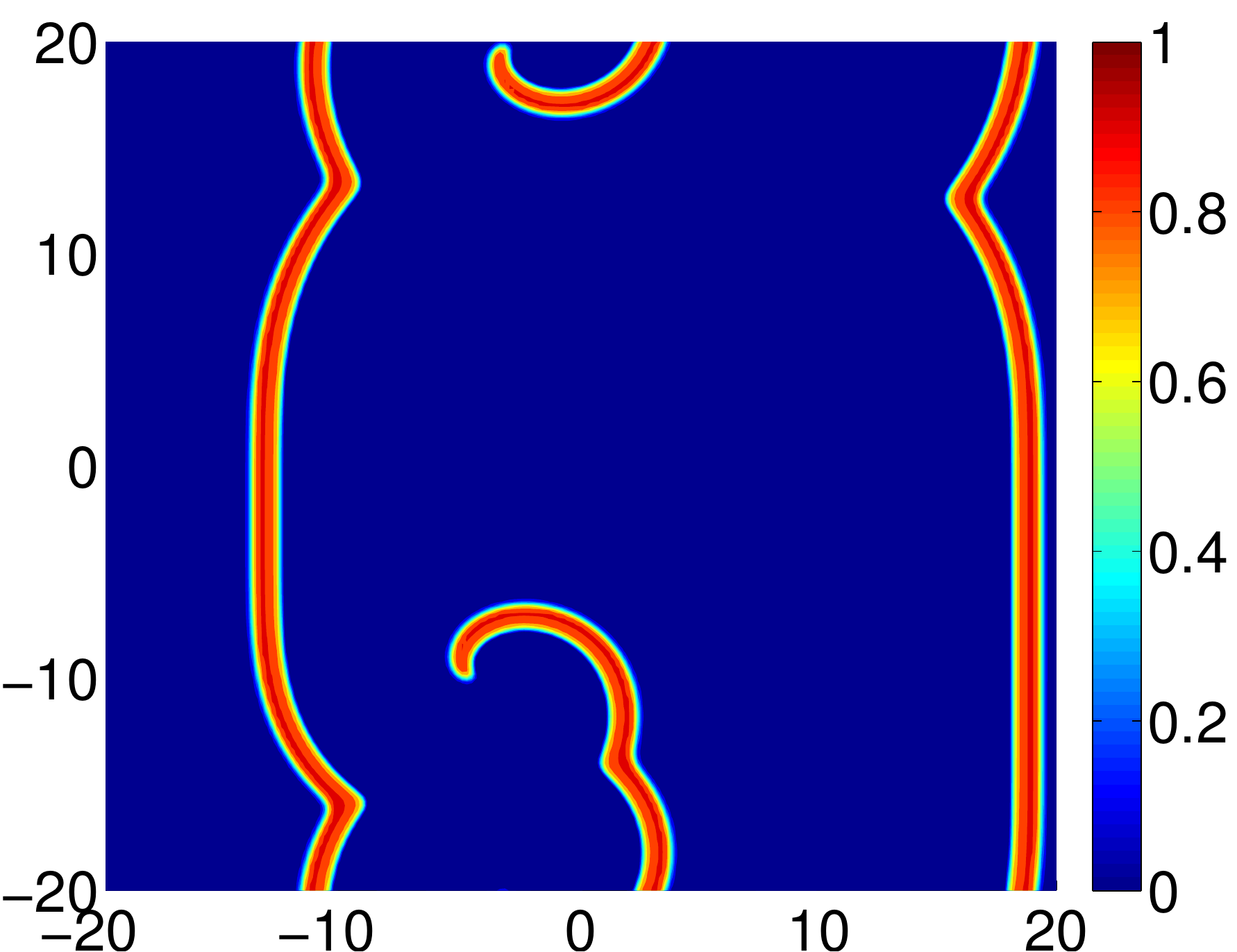}}
  %  \caption{Caption 1}
  (a2)
  \medskip
  \end{minipage}
  \hspace{0.03\linewidth}%
  \begin{minipage}[b]{0.3\textwidth}
    \centerline{
    \includegraphics[width=2.0in,angle=0,scale=1.0]{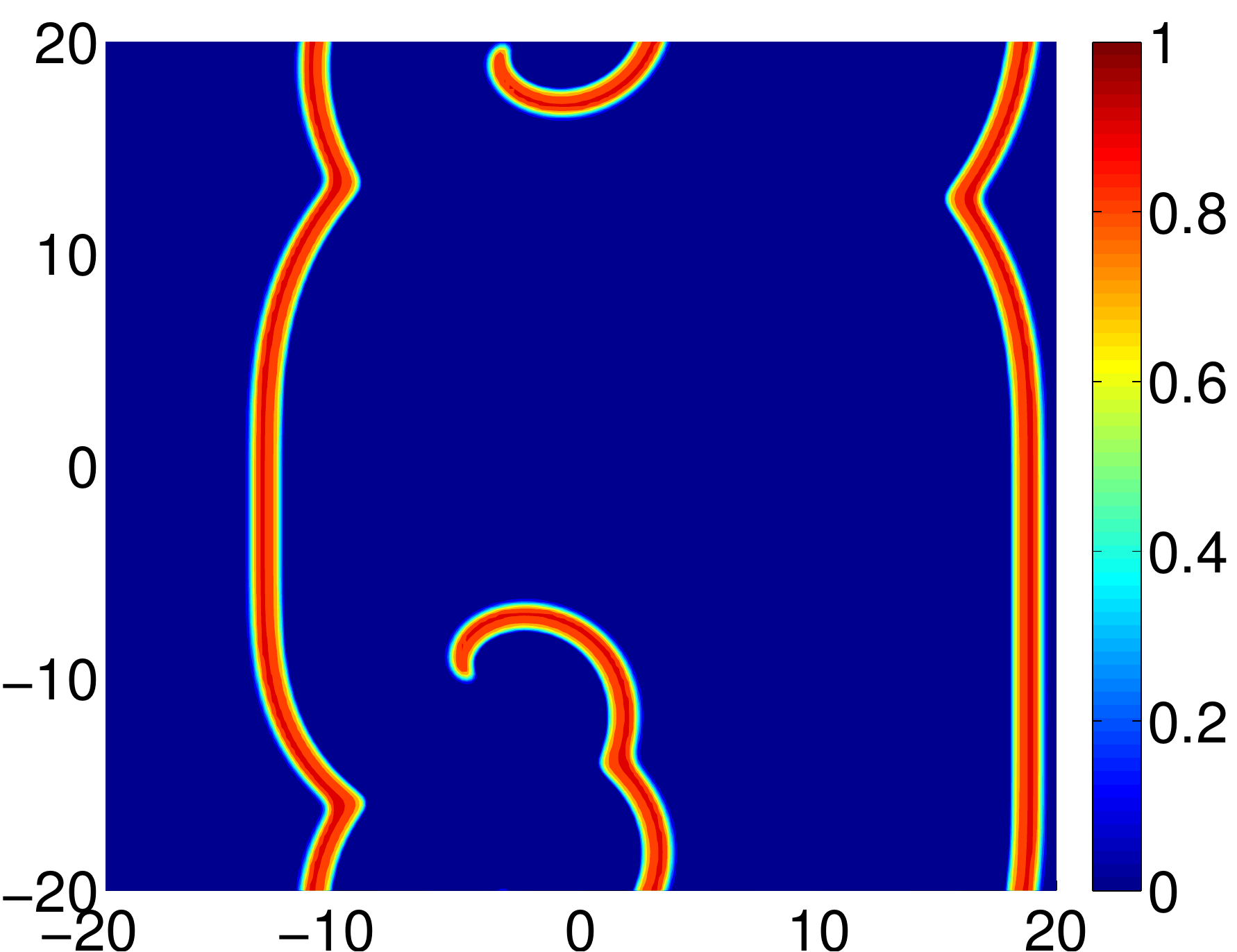}}
  %  \caption{Caption 1}
  (a3)
  \medskip
  \end{minipage}
  \\[0.5pt]%
  \hspace{0.1\textwidth}%
  \begin{minipage}[b]{0.3\textwidth}
    \centerline{
    \includegraphics[width=2.0in,angle=0,scale=1.0]{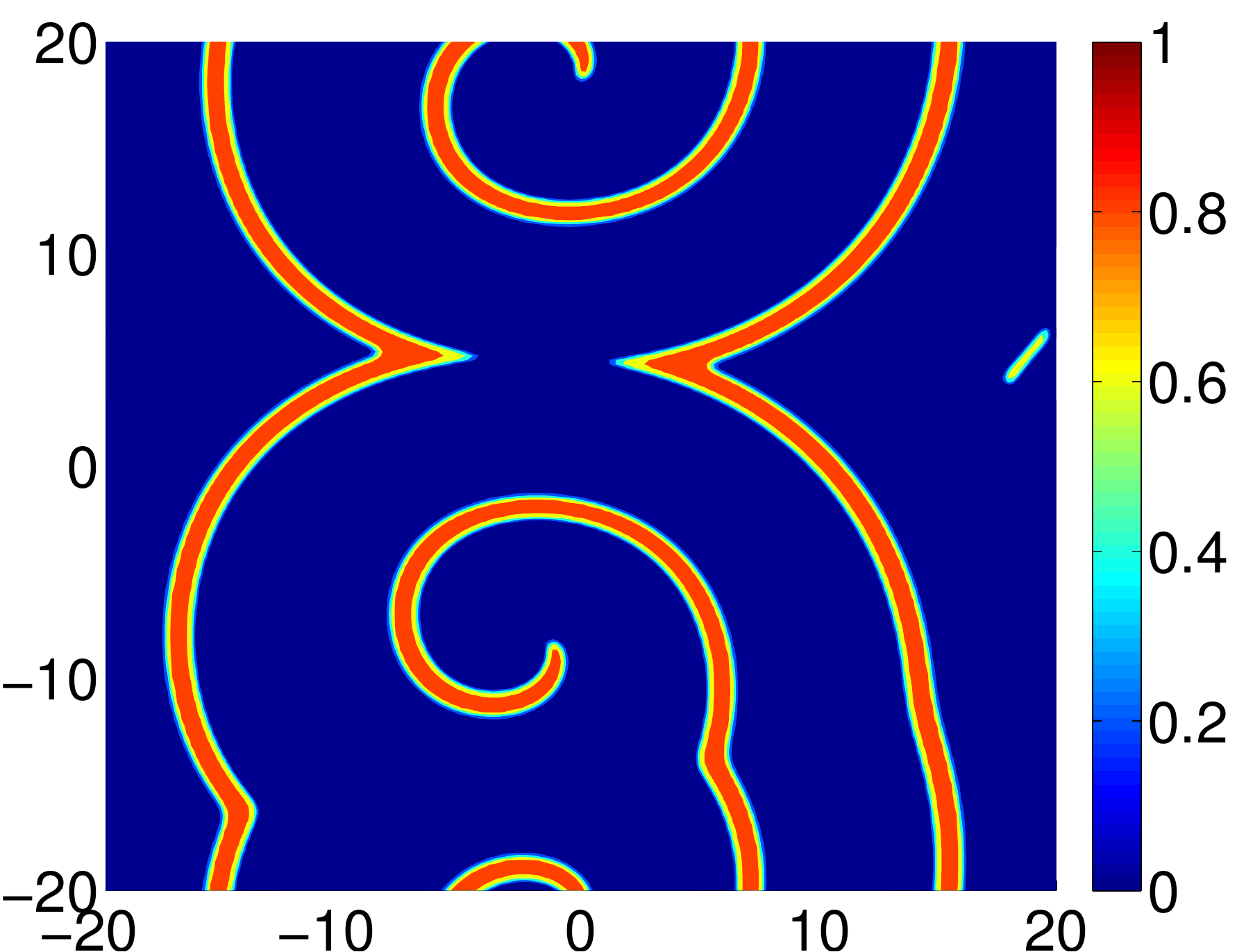}}
   % \caption{Caption 2}
   (b1)
   \medskip
  \end{minipage}
  \hspace{0.02\textwidth}
  \begin{minipage}[b]{0.3\textwidth}
    \centerline{
    \includegraphics[width=2.0in,angle=0,scale=1.0]{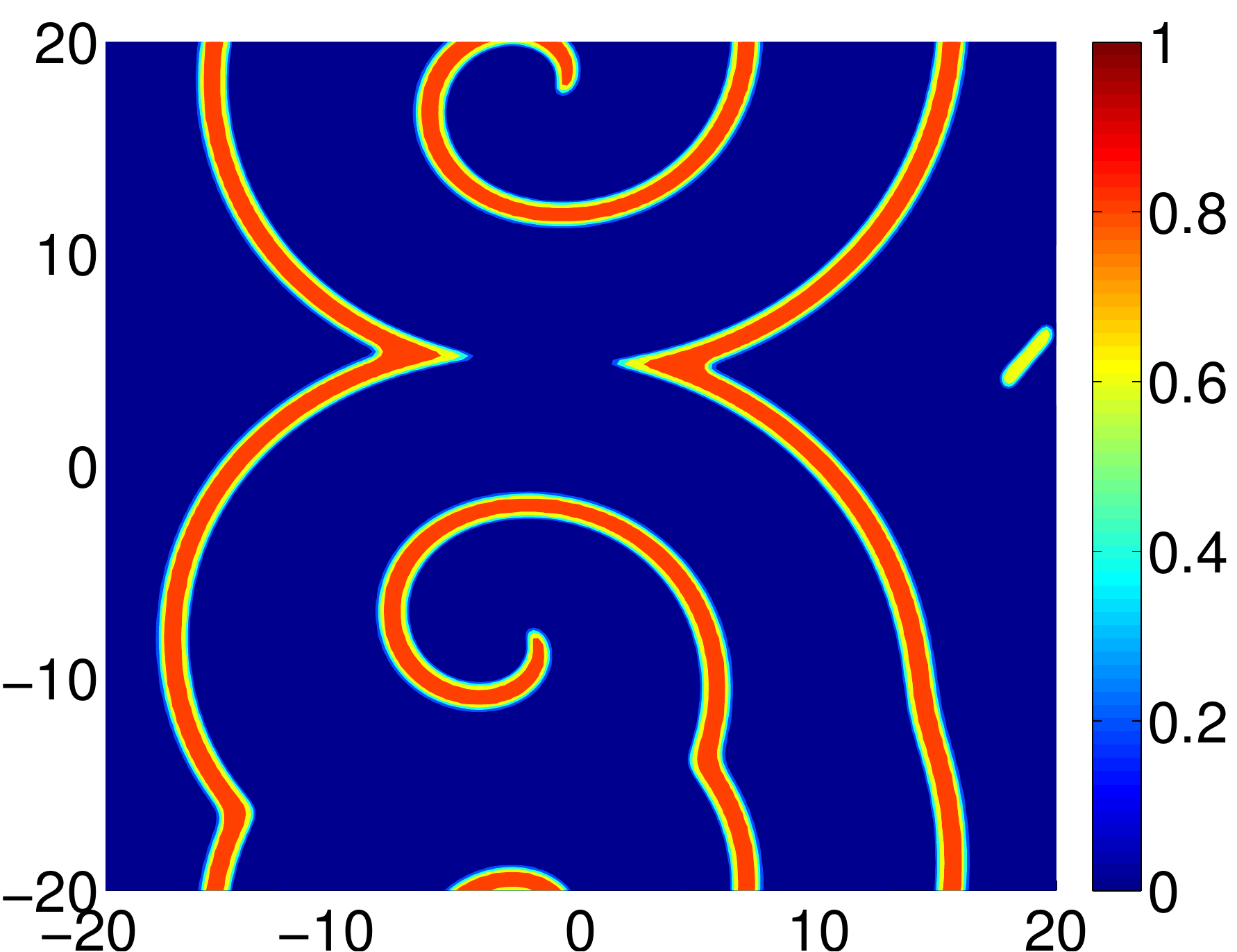}}
    %\caption{Caption 3}
    (b2)
    \medskip
  \end{minipage}
   \hspace{0.02\textwidth}
  \begin{minipage}[b]{0.3\textwidth}
    \centerline{
    \includegraphics[width=2.0in,angle=0,scale=1.0]{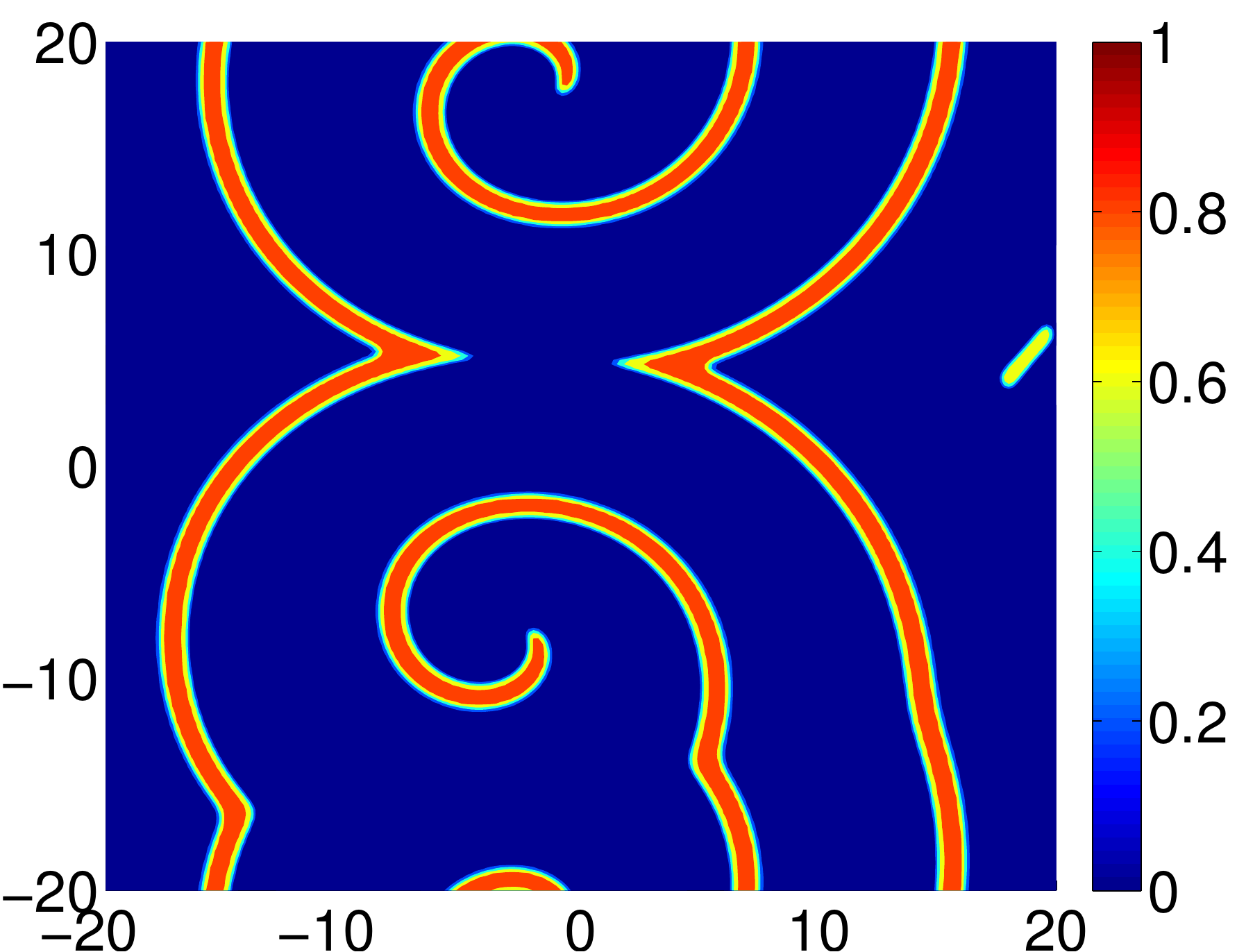}}
    %\caption{Caption 3}
    (b3)
    \medskip
  \end{minipage}
  \\[0.5pt]%
  \hspace{0.1\linewidth}%
  \begin{minipage}[b]{0.3\textwidth}
    \centerline{
    \includegraphics[width=2.0in,angle=0,scale=1.0]{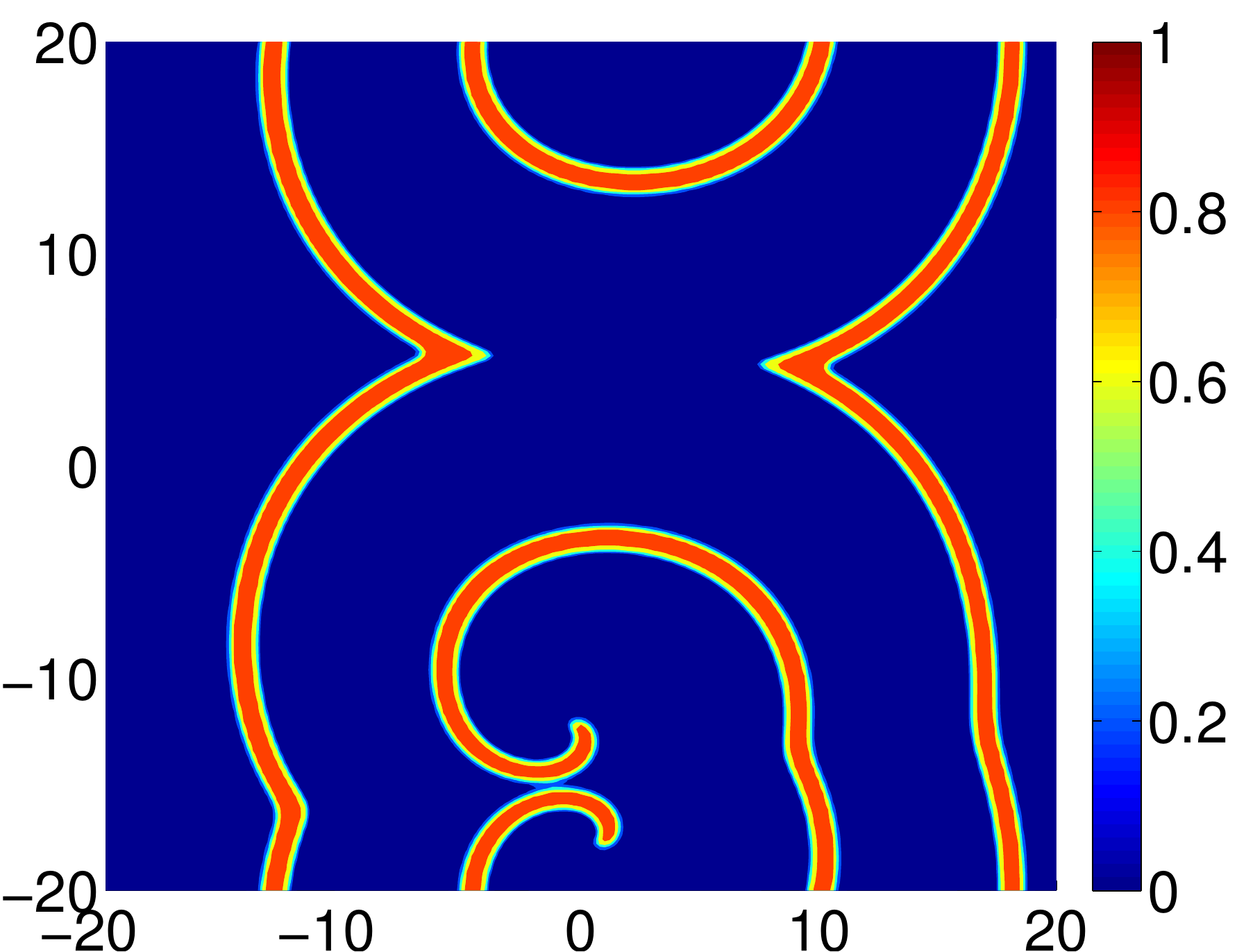}}
    %\caption{Caption 4}
    (c1)
    \medskip
  \end{minipage}%
  \hspace{0.03\linewidth}%
  \begin{minipage}[b]{0.3\textwidth}
    \centerline{
    \includegraphics[width=2.0in,angle=0,scale=1.0]{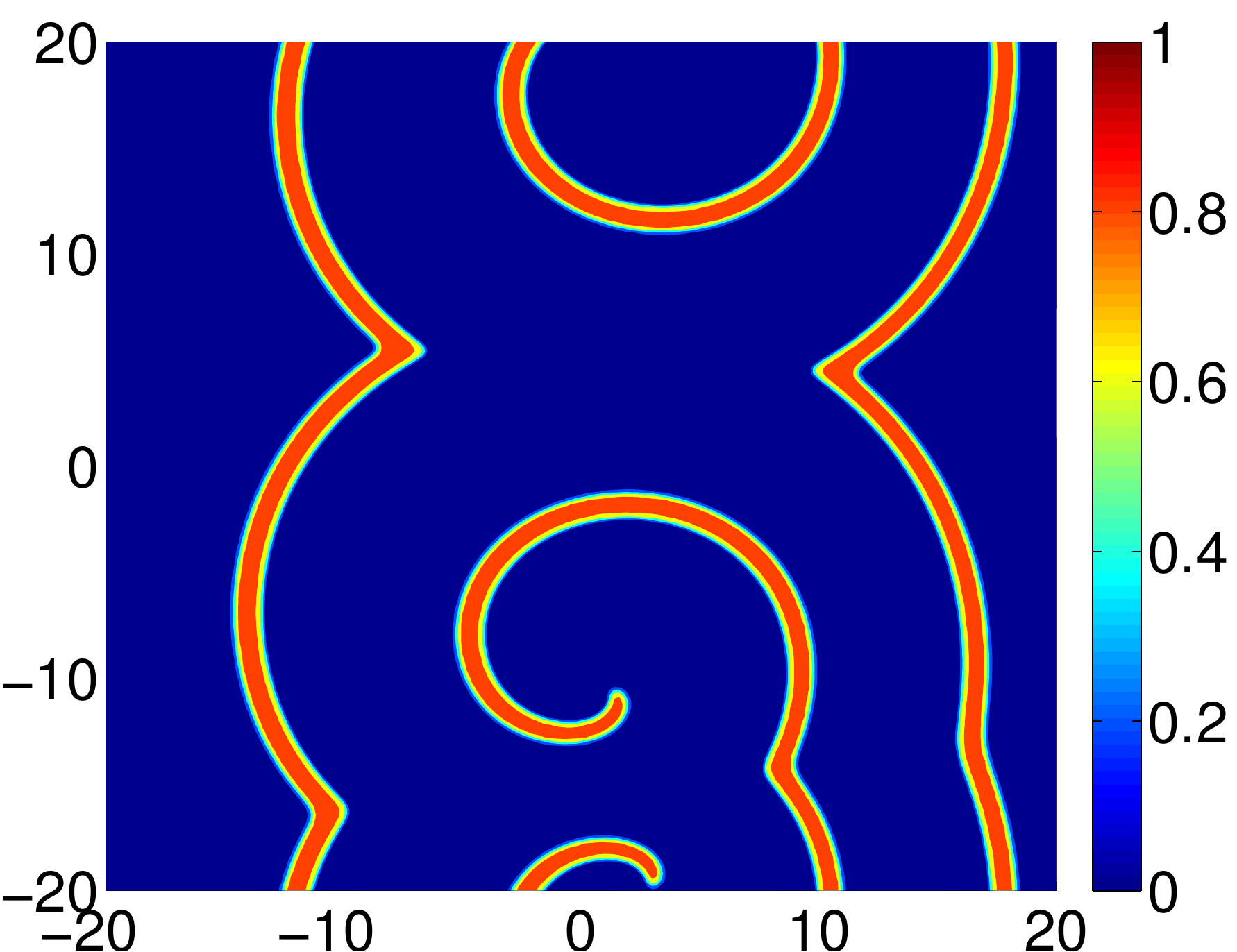}}
    % \caption{Caption 5}
    (c2)
    \medskip
  \end{minipage}
  \hspace{0.03\linewidth}%
  \begin{minipage}[b]{0.3\textwidth}
    \centerline{
    \includegraphics[width=2.0in,angle=0,scale=1.0]{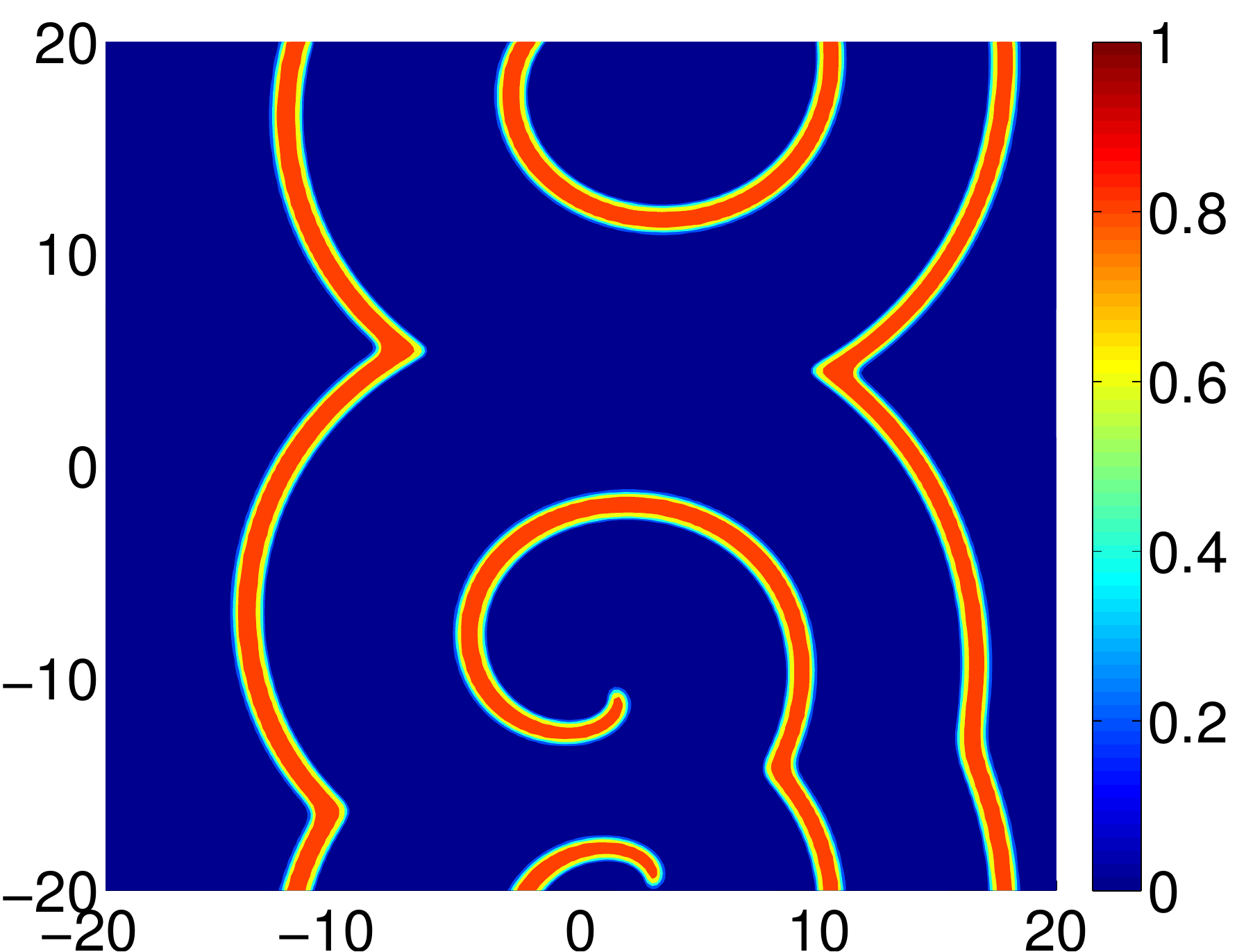}}
    % \caption{Caption 5}
    (c3)
    \medskip
  \end{minipage}
    \caption{Numerical simulations of the concentration of activator $u$ for Fitzhugh-Nagumo reaction-diffusion model at different times. (a1-a3) t = 2; (b1-b3) t = 5; (c1-c3) t = 10. (a1-c1) Lie-Trotter without corrector; (a2-c2) Lie-Trotter with two correctors; (a3-c3) Lie-Trotter with three correctors. } \label{fig_FN}
\end{figure}

%\pagebreak
\noindent {\bf Example 4.} \textbf{Fitzhugh-Nagumo reaction-diffusion model.}
A simple mathematical model of an excitable medium is Fitzhugh-Nagumo (FHN) equations \cite{Fife}. FHN equations with diffusion can be written as
\begin{align}
 \begin{cases}\label{FNeqn1}
 \vspace{0.1in}
 \displaystyle{\frac{\partial u}{\partial t}}=  D_u \nabla^2 u +\frac{1}{\delta}h(u,v),  \\
                  %               \vspace{0.1in}
 \displaystyle{\frac{\partial v}{\partial t}}=  D_v \nabla^2 v+g(u,v),
  \end{cases}
  \end{align}
where $D_u$, $D_v$ are the diffusion coefficients for activator $u$ and inhibitor $v$ respectively, and $\delta$ is a real parameter. We consider the classical cubic FHN local dynamics \cite{Keener, Olmos}
\begin{align}
 \begin{cases}\label{FNeqn2}
 \vspace{0.1in}
 \displaystyle   h(u,v) = Cu(1-u)(u-a)-v,  \\
                  %               \vspace{0.1in}
 \displaystyle  g(u,v) = u - d v,
  \end{cases}
  \end{align}
where $C$, $a$ and $d$ are dimensionless parameters. We perform the numerical experiment for (\ref{FNeqn1}) and (\ref{FNeqn2}) on the domain $[-20, 20]\times [-20, 20]$ with periodic boundary conditions. The parameters are chosen as following, $D_u = 1$, $D_v = 0$, $a = 0.1$, $C = 1$, $d = 0.5$, and $\delta = 0.005$. The initial condition is

\begin{equation}
   \label{FNeqn3}
\displaystyle u(x,y,0) = \left\{
  \begin{array}{ll}
  \vspace{0.05in}
     0, & \hbox{if $ \{x < 0 \}  \bigcup \{ y > 5 \}$;} \\
  \displaystyle  \frac{1}{(1+e^{4(|x|-5)})^2}- \frac{1}{(1+e^{4(|x|-1)})^2}, & \hbox{otherwise.}
  \end{array}
\right.
\end{equation}

\begin{equation}
   \label{FNeqn4}
\displaystyle v(x,y,0) = \left\{
  \begin{array}{ll}
  \vspace{0.05in}
     0.15, & \hbox{if $ \{x < 1 \}  \bigcap \{ y > -10 \}$;} \\
  0  , & \hbox{otherwise.}
  \end{array}
\right.
\end{equation}
 The domain is partitioned with a $200 \times 200$ grid. Figure \ref{fig_FN} shows the numerical solution to the concentration of the activator $u $ solving by Lie-Trotter splitting scheme, with three operators similar as Example 3. We observe the spiral waves at $T = 2, 5, 10$, which show a good agreement with the reference solutions. The computational step size is $\Delta t = 0.005$ in all cases. We remark that, because the lower order schemes suffer more from the numerical error of diffusion than higher order ones, the pattern for $T = 10$ look more consistent if we take a smaller computational time step size or a more refined mesh. Similar patterns can also be obtained by IDC-OS scheme based on Strang splitting.

\begin{figure}[!htb]
    \centerline{
    \includegraphics[width=2.0in,angle=0,scale=1.6]{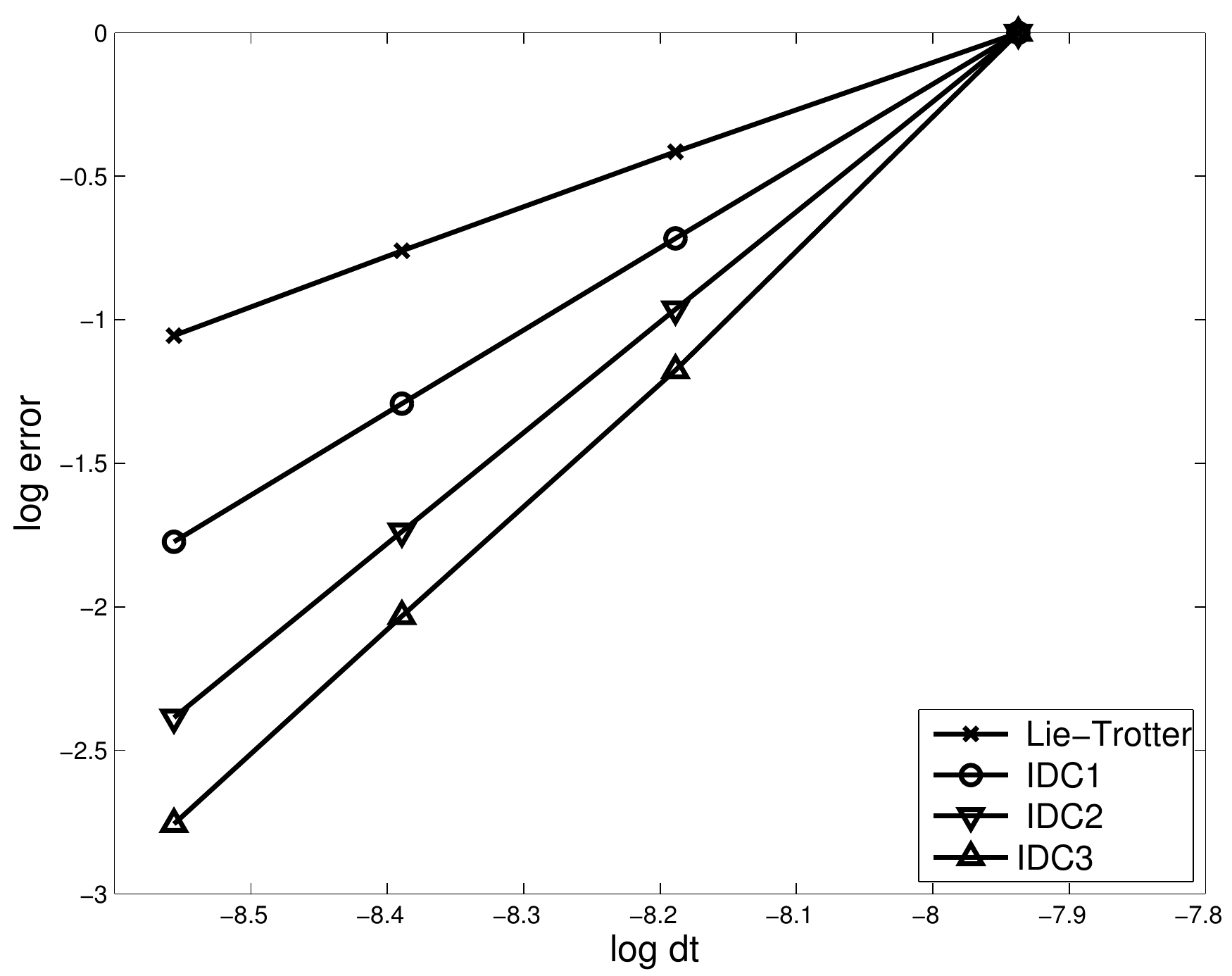}}
      \caption{Accuracy study for  Fitzhugh-Nagumo reaction-diffusion model. $t = 0.025$  } \label{fig_FN_order}
\end{figure}

 \noindent {\bf Remark.} Figure \ref{fig_FN_order}  shows the order of accuracy for Fitzhugh-Nagumo reaction-diffusion model (\ref{FNeqn1})  at $t = 0.025$, in which we clearly observe the order increase of IDC-OS scheme with successive correction steps. However, for a more stiff parameter such as $\delta = 10^{-10}$, order reduction phenomena is observed in the convergence study.  Similar observation is made in the Example 5 on Schnakenberg model.  How to approximate the residual integral and design a robust solver for stiff ODEs is an open question.  Recently work in \cite{Rokhlin}, the authors proposed a highly accurate solver based on an approximation of the integral of the residual  as a linear combination of exponentials on uniform quadrature nodes.  Their method is shown to do a good job of attaining high order of accuracy with correction steps and preserving the stability region of the original implicit time integrator which is used as the base scheme.

\begin{figure}[!htb]
  %\centering
  \begin{minipage}[b]{0.3\textwidth}
    \centerline{
    \includegraphics[width=2.0in,angle=0,scale=1.0]{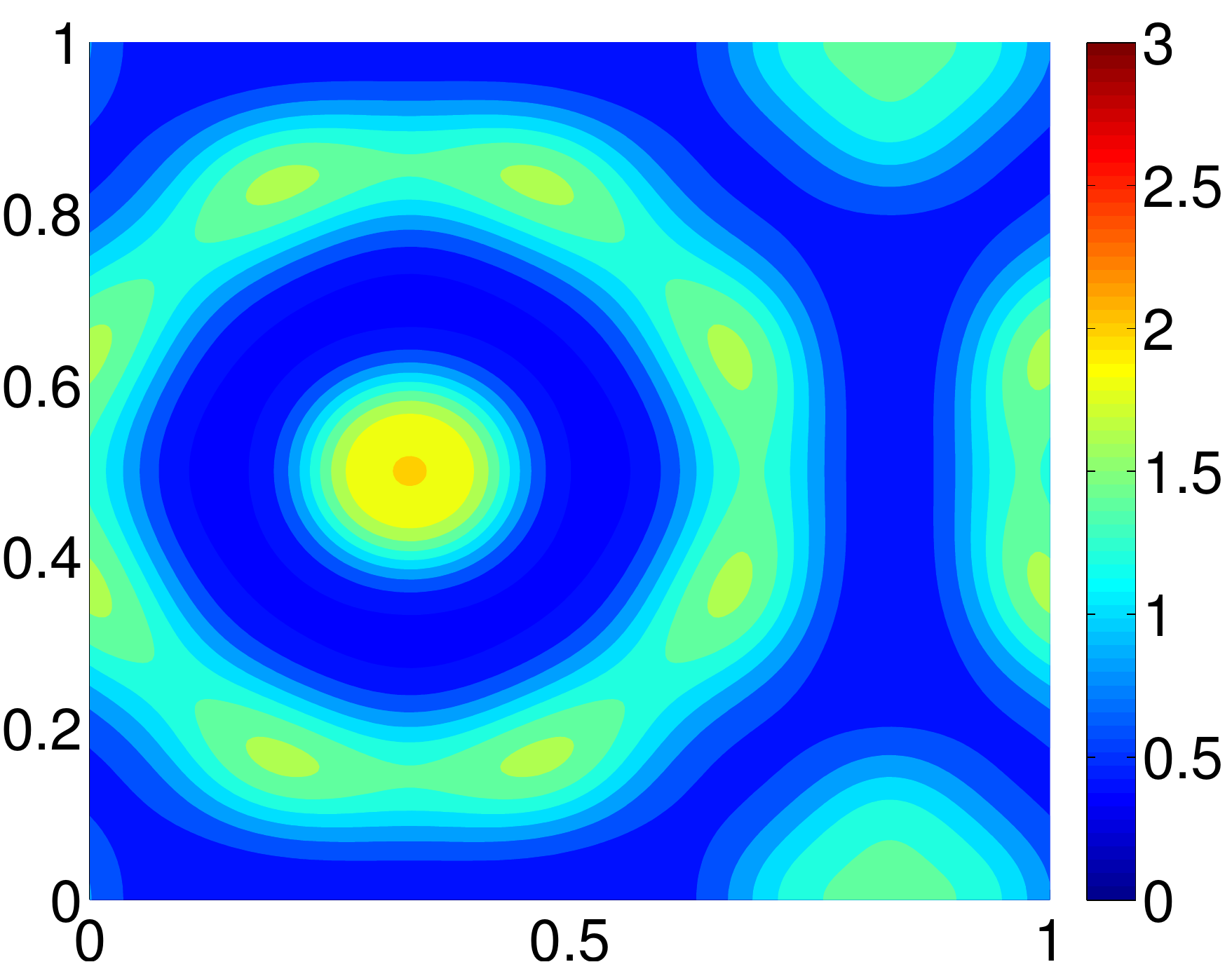}}
   % \caption{Caption 3}
      (a1)
    \medskip
  \end{minipage}%
  \hspace{0.04\linewidth}%
  \begin{minipage}[b]{0.3\textwidth}
    \centerline{
    \includegraphics[width=2.0in,angle=0,scale=1.0]{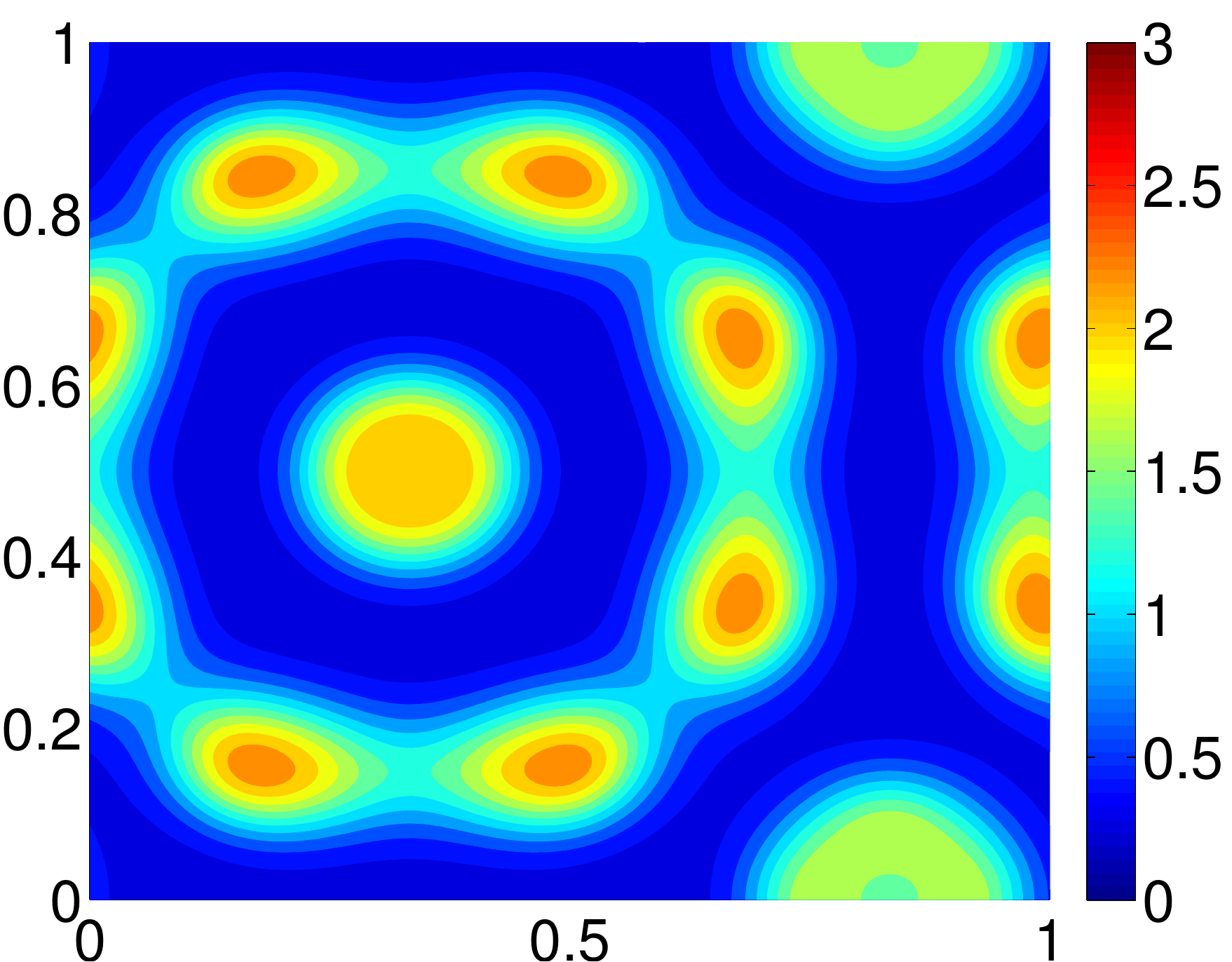}}
  %  \caption{Caption 1}
  (a2)
  \medskip
  \end{minipage}
  \hspace{0.03\linewidth}%
  \begin{minipage}[b]{0.3\textwidth}
    \centerline{
    \includegraphics[width=2.0in,angle=0,scale=1.0]{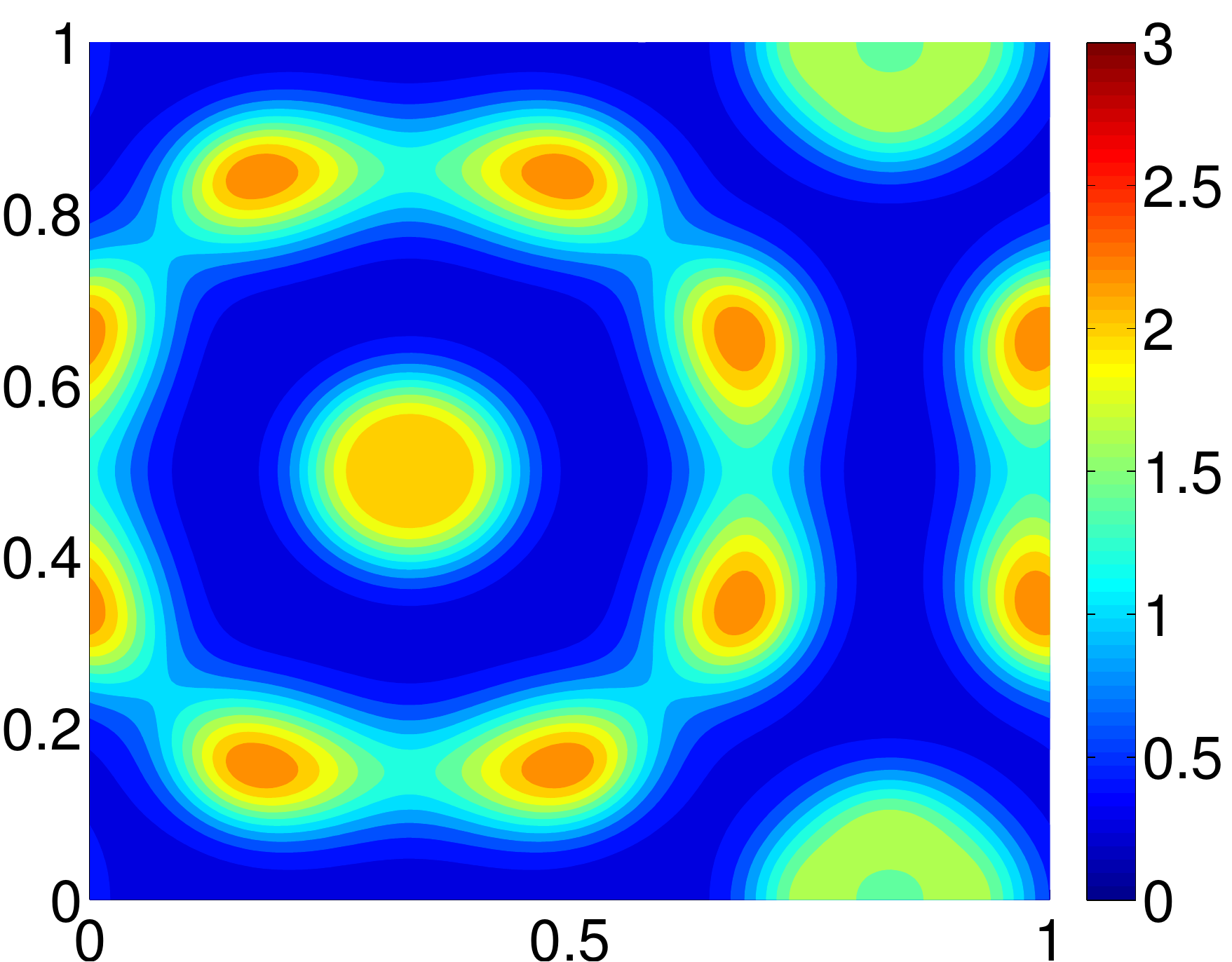}}
  %  \caption{Caption 1}
  (a3)
  \medskip
  \end{minipage}
  \\[0.5pt]%
  \hspace{0.1\textwidth}%
  \begin{minipage}[b]{0.3\textwidth}
    \centerline{
    \includegraphics[width=2.0in,angle=0,scale=1.0]{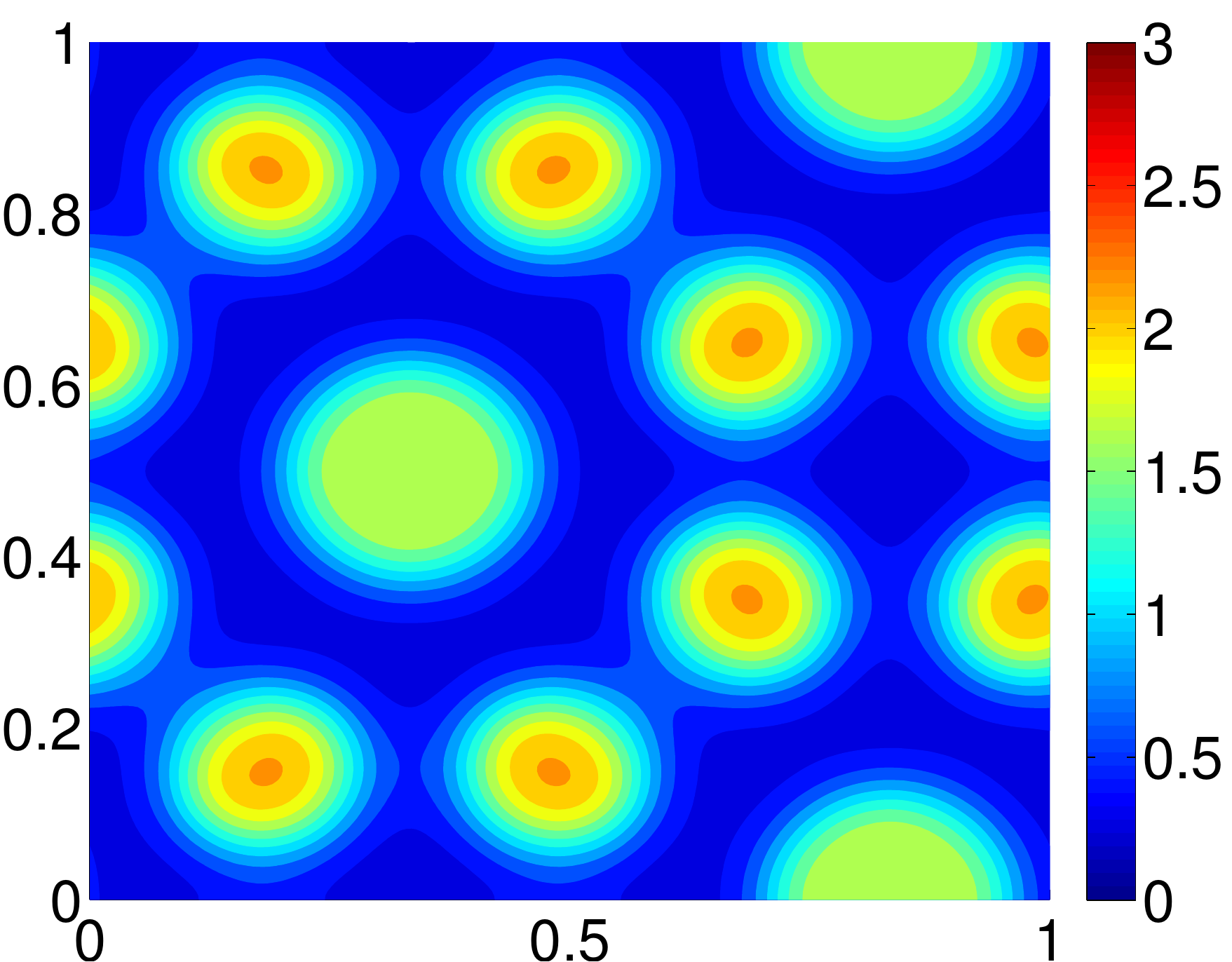}}
   % \caption{Caption 2}
   (b1)
   \medskip
  \end{minipage}
  \hspace{0.02\textwidth}
  \begin{minipage}[b]{0.3\textwidth}
    \centerline{
    \includegraphics[width=2.0in,angle=0,scale=1.0]{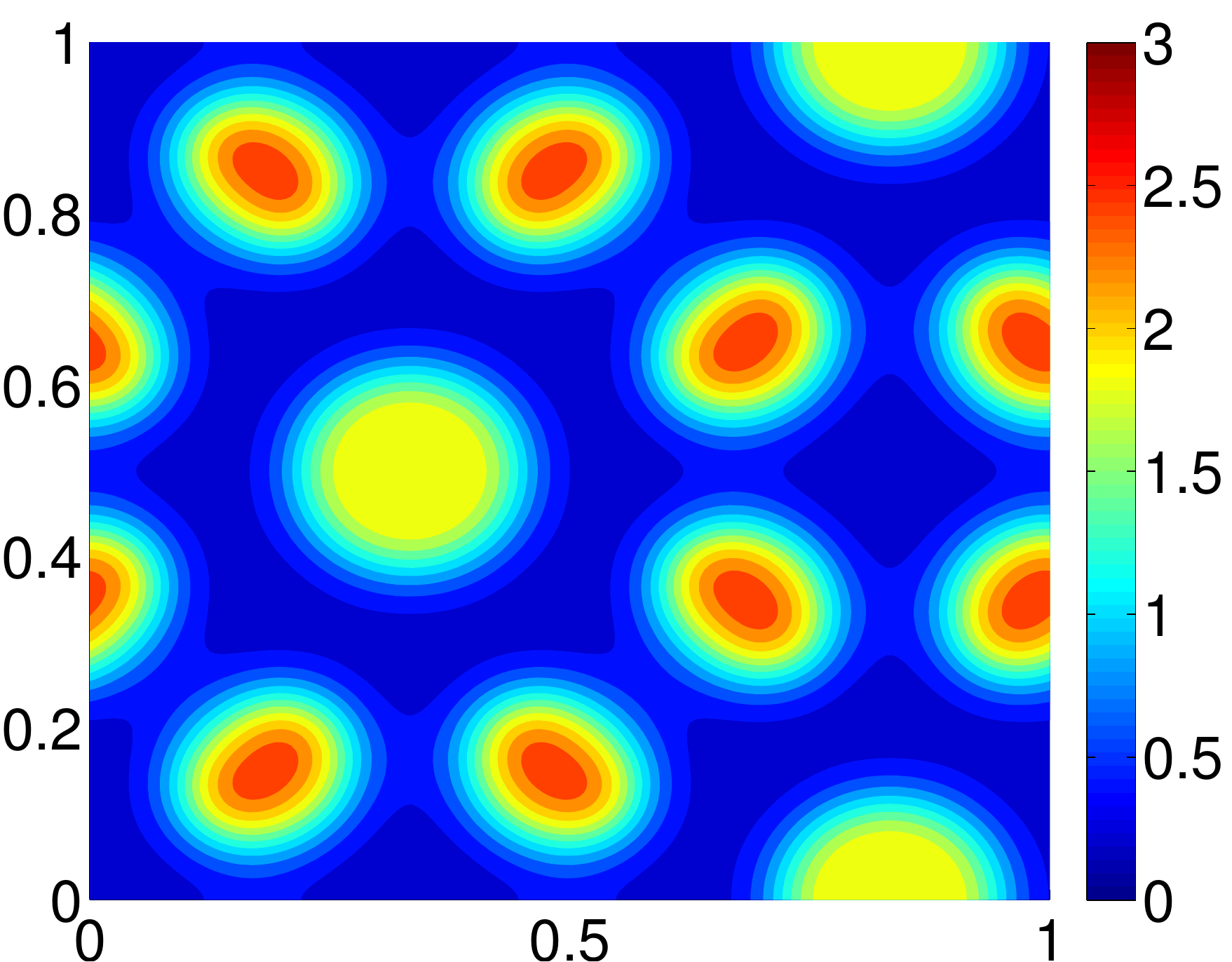}}
    %\caption{Caption 3}
    (b2)
    \medskip
  \end{minipage}
   \hspace{0.02\textwidth}
  \begin{minipage}[b]{0.3\textwidth}
    \centerline{
    \includegraphics[width=2.0in,angle=0,scale=1.0]{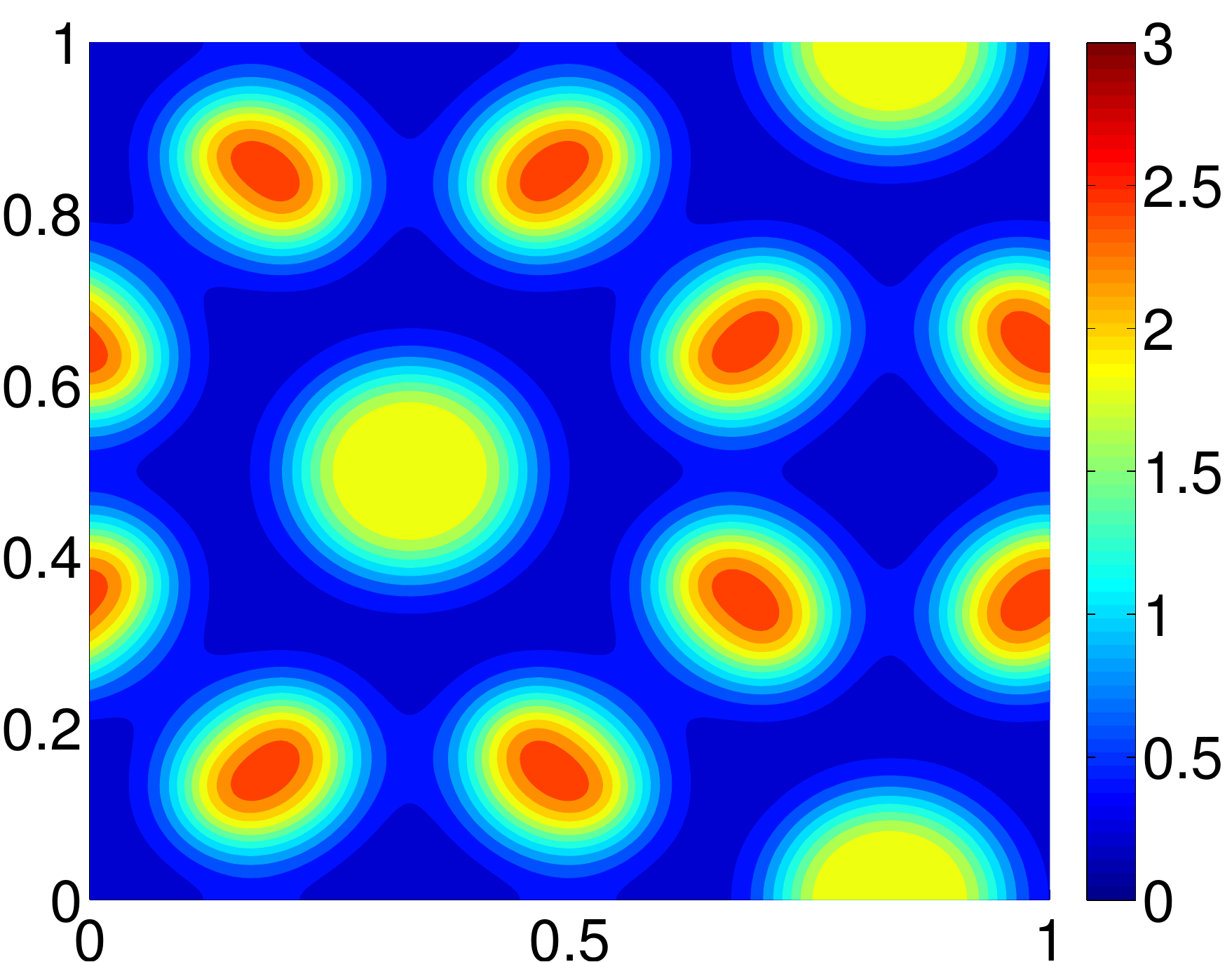}}
    %\caption{Caption 3}
    (b3)
    \medskip
  \end{minipage}
  \\[0.5pt]%
  \hspace{0.1\linewidth}%
  \begin{minipage}[b]{0.3\textwidth}
    \centerline{
    \includegraphics[width=2.0in,angle=0,scale=1.0]{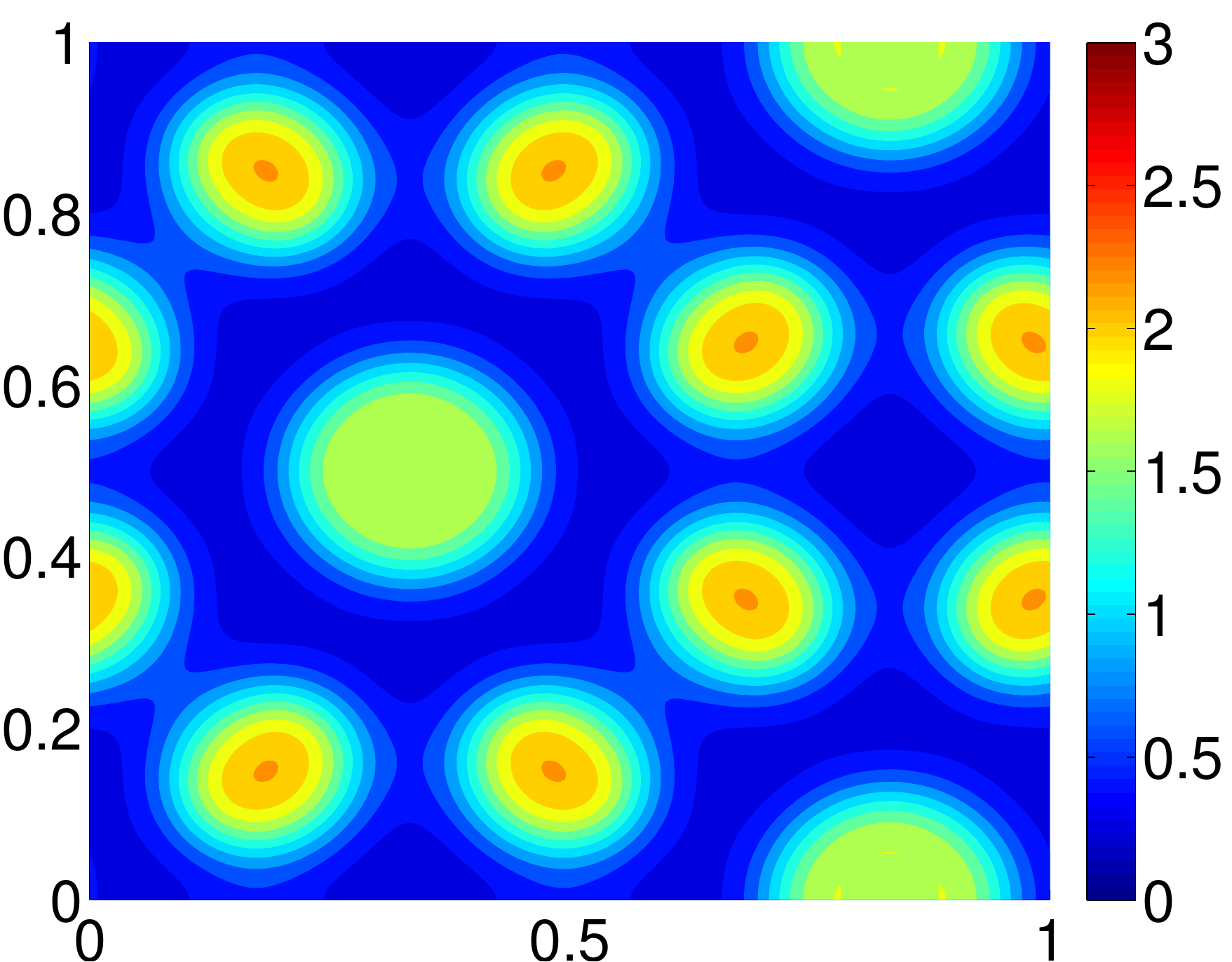}}
    %\caption{Caption 4}
    (c1)
    \medskip
  \end{minipage}%
  \hspace{0.04\linewidth}%
  \begin{minipage}[b]{0.3\textwidth}
    \centerline{
    \includegraphics[width=2.0in,angle=0,scale=1.0]{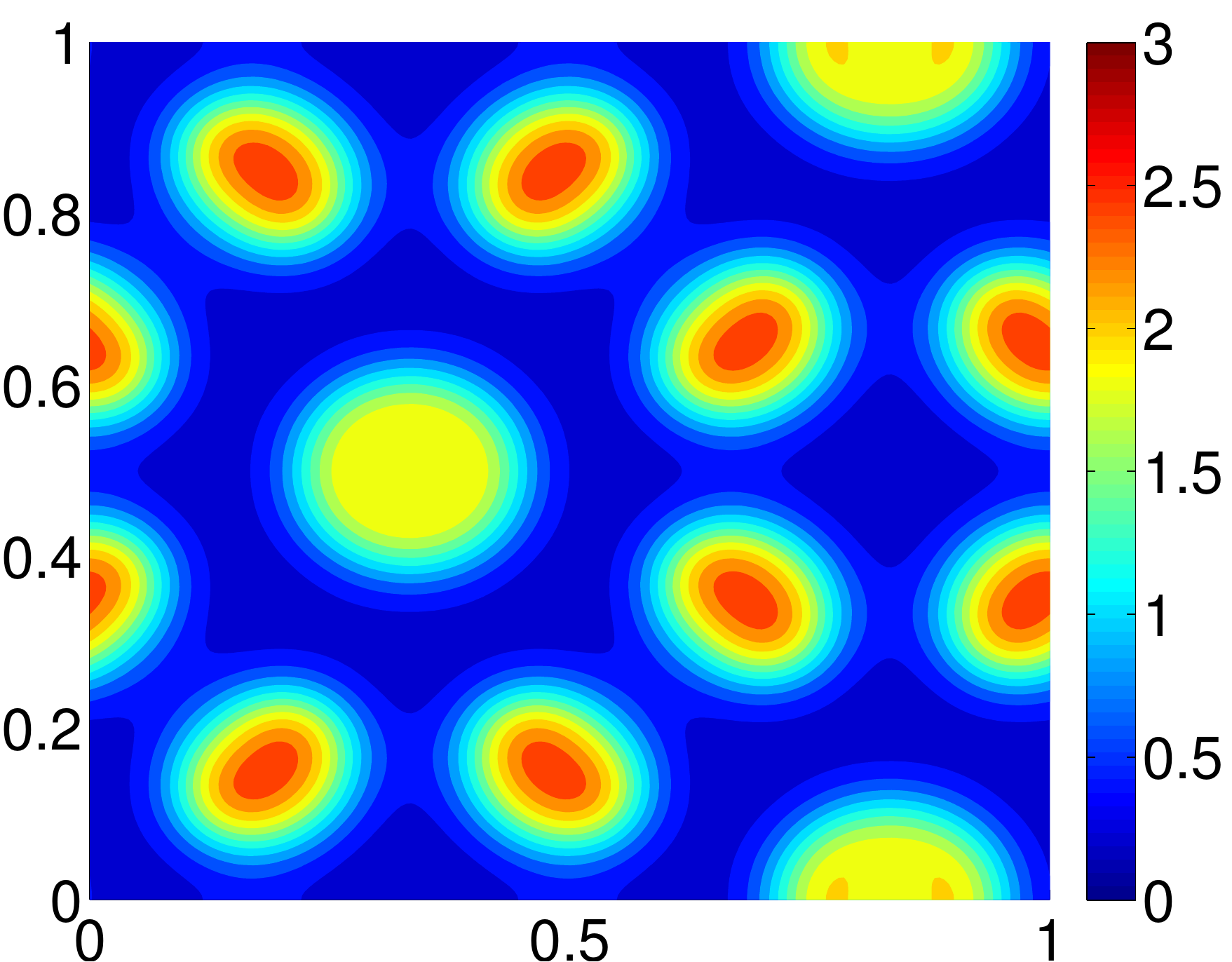}}
    % \caption{Caption 5}
    (c2)
    \medskip
  \end{minipage}
  \hspace{0.03\linewidth}%
  \begin{minipage}[b]{0.3\textwidth}
    \centerline{
    \includegraphics[width=2.0in,angle=0,scale=1.0]{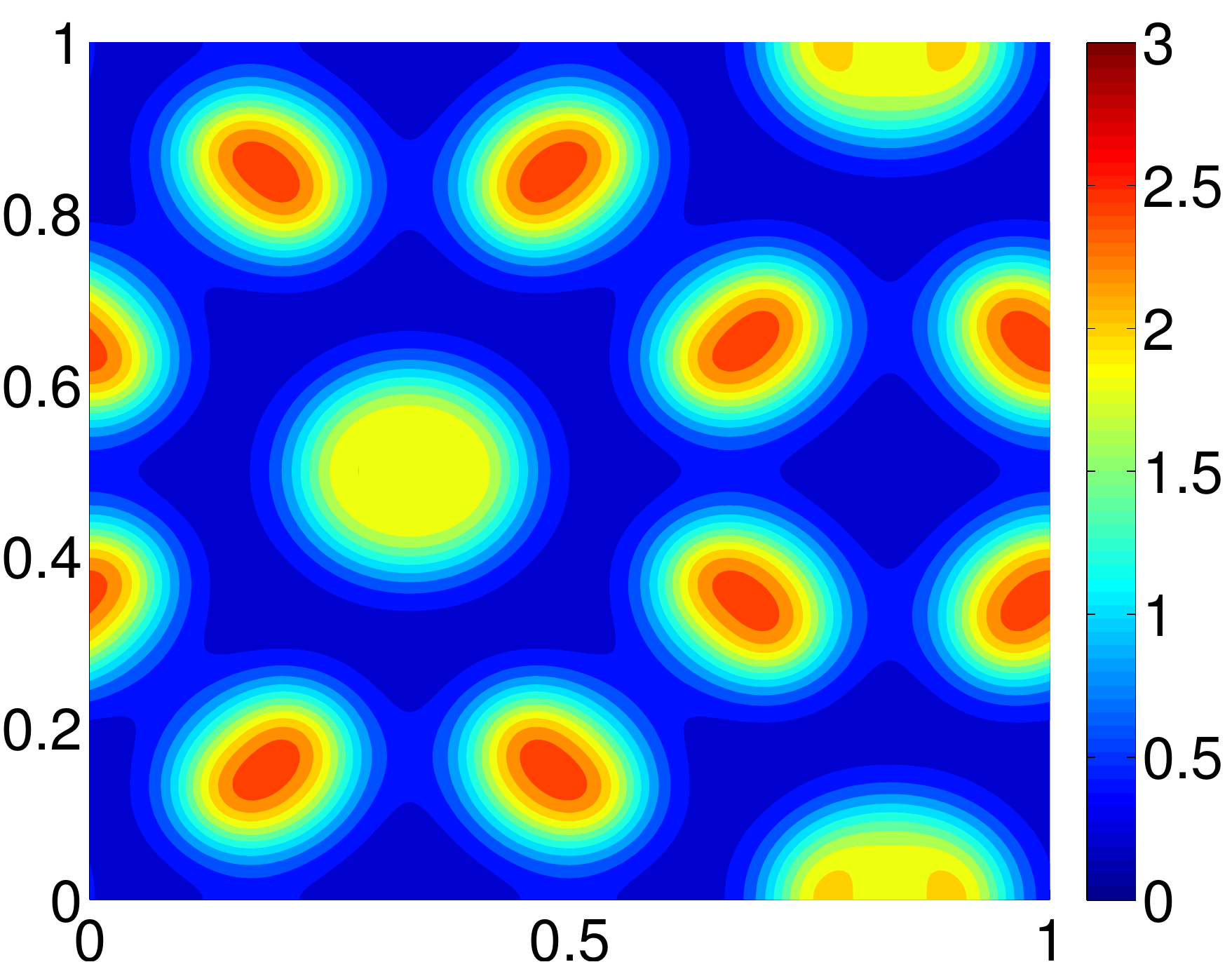}}
    % \caption{Caption 5}
    (c3)
    \medskip
  \end{minipage}
    \caption{Numerical simulations of the concentration of activator $C_a$ for Schnakenberg reaction-diffusion model at different times. (a1-a3) t = 0.5; (b1-b3) t = 1; (c1-c3) t = 1.5. (a1-c1) Lie-Trotter without corrector; (a2-c2) Lie-Trotter with one corrector; (a3-c3) Lie-Trotter with two correctors. } \label{fig_pig}
\end{figure}

%%%%%%%%%%%%%%%%%%%%%%%%%%%%%%%%%%%%%%%%%%%%%%%%%%%%%%%%%%%%%%
%\pagebreak
%\newpage
\noindent {\bf Example 5.} \textbf{Schnakenberg model.} The Schnakenberg system \cite{Schnakenberg} has been used to model the spatial distribution of a morphogen. It has the following form
\begin{align}
 \begin{cases}\label{SMeqn1}
 \vspace{0.1in}
 \displaystyle{\frac{\partial C_a}{\partial t}}=  D_1 \nabla^2 C_a+\kappa(a-C_a+C_a^2C_i),  \\
                  %               \vspace{0.1in}
 \displaystyle{\frac{\partial C_i}{\partial t}}=  D_2 \nabla^2 C_i+\kappa(b-C_a^2C_i),
  \end{cases}
  \end{align}
where $C_a$ and $C_i$ represents the concentration of activator and inhibitor, with $D_1$ and $D_a$ as the diffusion coefficients respectively. $\kappa$, $a$ and $b$ are rate constants of biochemical reactions. Following the setup in \cite{Hundsdorfer}, we take the initial conditions as
\begin{align}\label{SMeqn2}
\displaystyle C_a(x,y,0)& = a + b + 10^{-3}e^{-100((x-\frac{1}{3})^2+(y-\frac{1}{2})^2)}, \\
\displaystyle C_i (x,y,0)& =  \frac{b}{(a+b)^2},
\end{align}
and the boundary conditions are periodic. The parameters are $\kappa = 100$, $a = 0.1305$, $b = 0.7695$, $D_1 = 0.05$ and $D_2 = 1$. The computational domain is $[0, 1] \times [0, 1]$. The numerical simulations with Lie-Trotter splitting is performed on a $200 \times 200$ spatial grid and the numerical dynamical process of the concentration of the activator $C_a$ at different times are shown in Figure \ref{fig_pig}, we can observe that the initial data are amplified and spreads, leading to thew formation of spot pattern. The computational time step size is chosen as $\Delta t = 0.001$. We also note that similar patterns can be obtained by IDC-OS scheme based on Strang splitting.

%%%%%%%%%%%%%%%%%%%%%%%%%%%%%%%%%%%%%%%%%%%%%%%%%%%%%%%%%%%%%%%
%%\pagebreak
%%\newpage
%\noindent {\bf Example 6.} \textbf{Van der Pol equation.}  In this example, we examine the nonlinear oscillatory test problem, Van der Pol's equation, 
%\begin{align}
% \begin{cases}\label{Vanderpol}
% \vspace{0.1in}
% \displaystyle{\frac{\partial u}{\partial t}}=  v,  \\
%                  %               \vspace{0.1in}
% \displaystyle{\frac{\partial v}{\partial t}}=  \frac{1}{\epsilon}(-u+(1-u^2)v),
%  \end{cases}
%  \end{align}
%with initial condition as $u(0) = 2$, $v(0) = -\frac{2}{3}+\frac{10}{81}\epsilon -\frac{292}{2187}\epsilon^2-\frac{1814}{19683}\epsilon^3$ and end time is taken as $T = 0.5$. Consistent with previous examples, we propose
%to split the right-hand side functions evenly.
%
%

% section 7
%\input{conclusion.tex}

\section{Conclusion}
In this paper, we have provided a general temporal framework for the construction of high order operator splitting methods based on the integral deferred correction procedure. The method can achieve arbitrary high order via solving correction equation whereas reduce the computational cost by taking the advantage of operator splitting. Error analysis and  numerical examples for IDC-OS methods are performed to show that the proposed IDC framework successfully enhances the order of accuracy in time.  A study on order reduction for very stiff problems will be part of our future work.

% section 8
%\input{acknowledgment.tex}

\section*{Acknowledgments}
AJC supported in part by AFOSR grants FA9550-11-1-0281, FA9550-12-1-0343 and
FA9550-12-1-0455, NSF grant DMS-1115709, and MSU Foundation grant SPG-RG100059.  ZX is supported by NSF grant DMS-1316662. Additionally, the authors would like to thank Prof. William Hitchon and Dr. David Seal for helpful comments on this work.

\bibliographystyle{abbrv}
\bibliography{BigBib}
\end{document}